\documentclass[11pt]{article}
\usepackage{amsmath,amsfonts,amssymb,mathrsfs,cases,mathdots,colordvi,url}
\usepackage{color}
\newtheorem{defn}{Definition}[section]

\newtheorem{thm}{Theorem}[section]
\newtheorem{prop}{Proposition}[section]
\newtheorem{cor}{Corollary}[section]
\newtheorem{example}{Example}[section]
\newtheorem{remark}{Remark}[section]

\renewcommand{\Box}{\rule{2.2mm}{2.2mm}}

\newcommand{\K}{{\cal K}}

\newcommand{\R}{\mathbb{R}}
 \def\beginproof{\par\noindent {\bf Proof.}\ \ }
 \def\endproof{\hskip .5cm $\Box$ \vskip .5cm}

\def\beginproof{\par\noindent {\bf Proof.}\ \ }
\def\endproof{\hskip .5cm $\Box$ \vskip .5cm}
\newcommand{\ds}{\displaystyle}

\begin{document}
\title{Verifiable sufficient conditions for the error bound property of second-order cone complementarity problems}
\author{Jane J. Ye\thanks{Department of Mathematics
and Statistics, University of Victoria, Victoria, B.C., Canada V8W 2Y2, e-mail: janeye@uvic.ca. The research of this author was partially
supported by NSERC.}\  \ \ \ and \ \ Jinchuan Zhou\thanks{Department of Statistics, School of Mathematics
and Statistics, Shandong University of Technology,
 Zibo 255049, P.R. China, e-mail: jinchuanzhou@163.com. This author's work is supported by National
 Natural Science Foundation of China (11771255, 11101248) and Shandong Province Natural Science Foundation (ZR2016AM07).}}

\maketitle
\begin{abstract} The  error bound property for a solution set defined by a set-valued mapping refers to an inequality that bounds the distance between vectors closed to a solution of the given set by a residual function. The  error bound property is a Lipschitz-like/calmness property of the perturbed solution mapping, or equivalently the  metric subregularity of the  underlining set-valued mapping. It has been proved to be extremely useful in analyzing the convergence  of many algorithms for solving optimization problems, as well as serving as a constraint qualification for optimality conditions. In this paper, we study the  error bound property
for the solution  set of a very general second-order cone complementarity problem (SOCCP). We derive some sufficient conditions for error bounds for SOCCP which is verifiable based on the initial problem data.

\vskip 10 true pt

\vskip 10 true pt


\end{abstract}

\section{Introduction}

In this paper we consider a second-order cone complementarity  problem (SOCCP) of finding $z\in \mathbb{R}^n$ satisfying the second-order cone complementarity system defined as
\begin{eqnarray}
&&{\cal K}  \ni G(z)  \perp H(z)  \in {\cal K}, \label{system-1}\\
&&  F(z) \in \Lambda, \label{system-2}
\end{eqnarray}
 where % $g:\mathbb{R}^n\to \mathbb{R}^p$, $h:\mathbb{R}^n\to \mathbb{R}^q$,
$\Lambda$ is a closed subset of $\mathbb{R}^l$,  $F:\mathbb{R}^n\to \mathbb{R}^l$,
 $G:\mathbb{R}^n\to \mathbb{R}^{m}$, $H:\mathbb{R}^n\to \mathbb{R}^{m}$ are continuously differentiable, $ a \perp b$ means that  vector $a$ is perpendicular to vector $b$, ${\cal K}$ is the Cartesian product of finitely many second-order cones (also called Lorentz cones), i.e.,
 \[
 {\cal K}:={\cal K}_1 \times {\cal K}_2 \times \cdots  \times {\cal K}_J,
 \]
with ${\cal K}_i:= \{x=(x_1,x_2)\in \mathbb{R}\times \mathbb{R}^{m_i-1}|\ x_1\geq \|x_2\|\}$  being the $m_i$-dimensional second-order cone and  $m=\sum_{i=1}^J m_i$.
 %$\|\cdot\|$ denoting the Euclidean norm and
 %${\cal K}_i$ with $m_i=1$ denoting the set of  nonnegative reals $\mathbb{R}_+$.

 One of the sources of the second-order cone complementarity system is the Karush-Kuhn-Tucker (KKT) optimality condition  for the second-order cone programming (see e.g. \cite{AG03,BR05,CT}), and the other is the equilibrium system for a Nash game where the constraints involving second-order cones (see e.g. \cite{HYF}).

  Let ${\cal F}$ denote the solution set of an SOCCP which contains all $z$ satisfying the second-order cone complementarity system (\ref{system-1})-(\ref{system-2}). In this paper, we study the following error bound property.
  We say that the second-order cone complementarity system has a local error bound
at $z^*\in {\cal F}$ if there exist a constant $\kappa >0$ and $U$ a neighborhood of $z^*$ such that
 \begin{eqnarray}
 && d(z, {\cal F})\leq \kappa \left \{d_\Lambda(F(z))+ \sum_{i=1}^J d_{\Omega_i} (G_i(z),H_i(z))\right \} \qquad  \forall z\in U, \label{errorb}
 \end{eqnarray}
 where
 $\Omega_i:=\{(x,y)|{\cal K}_i  \ni x \perp y  \in {\cal K}_i\}$
 is the  $m_i$-dimensional second-order cone complementarity set.  The right hand side of the inequality (\ref{errorb})  is a residual function, and hence  the existence of a local error bound  enables us to use the residual to measure the distance from a point $z$ that  is sufficiently close to $z^*$  to the solution set ${\cal F}$.
 It is easy to verify that the error bound property at $z^*$ is equivalent to the calmness of the set-valued mapping defined by
$$
{\cal F}(\alpha, \beta,\gamma):=\left \{z \Big |  \begin{array}{l}
{\cal K}  \ni (G(z)+\alpha)  \perp (H(z)+\beta)  \in {\cal K}\\
F(z)+\gamma\in \Lambda
\end{array} \right \}
$$ at ${(0,0,0,z^*)}\in {\rm gph}{\cal F}$. Since ${\cal F}(0, 0,0)={\cal F}$, the solution to  the second-order cone complementarity system, the set-valued map ${\cal F}(\alpha, \beta,\gamma)$ is the perturbed  solution mapping.  Hence, the calmness property is a Lipschitz-like property of the perturbed  solution mapping: there exist a constant $\kappa >0$,  $U$ a neighborhood of $z^*$, $W$ a neighborhood of $(0,0,0)$  such that
  $${\cal F}(\alpha,\beta,\gamma)
  \cap U\subseteq
   {\cal F} +
   \kappa
   \|(\alpha,\beta,\gamma)\| \bar{B} \quad \forall (\alpha,\beta,\gamma)\in W. $$

{For $i=1,\dots,J$}, it is easy to verify that
\begin{equation}
(x,y)\in \Omega_i \Longleftrightarrow x=\Pi_{{\cal K}_i}(x-y),\label{relation-1}
\end{equation}
and the following inequality holds:
$$d_{\Omega_i} (x,y)\leq \sqrt{2} \|x-\Pi_{{\cal K}_i}(x-y)\|, \quad \forall x,y\in \mathbb{R}^{m_i}, $$
where $\Pi_{{\cal K}_i}(z)$ denotes the metric projection  of $z$ onto ${\cal K}_i$. Therefore,  if the error bound property holds with the residual function $d_\Lambda(F(z))+ \sum_{i=1}^J d_{\Omega_i} (G_i(z),H_i(z))$, then the error bound property also holds with the natural residual function $d_\Lambda(F(z))+ \sum_{i=1}^J \|G_i(z)-\Pi_{{\cal K}_i}(G_i(z)-H_i(z))\| $.

The  error bound property and equivalently the calmness property is a very  important property. One of the applications of such a property is  the analysis of certain algorithms for solving the second-order cone complementarity  problem. In particular, it has recently been discovered that the condition that is crucial to the quadratic convergence of the Newton-type method is not the nonsingularity of the Jacobian {\it per se}, but rather one of its consequences--the error bound
%t is an essential condition to gurantee the quadratic convergence of a Newton-type method for the solution of a nonsmooth equation
property; see \cite{FFH}. Another application is the constraint qualification for the mathematical program with second-order cone complementarity constraints (SOC-MPCC); see \cite{YZ15}.

Although the  error bound property is an important property, there are very few results on sufficient conditions for the existence of error bounds, and these results are abstract  and not easy to verify; see e.g. \cite{FHKO,NZ00,Pang97,WuYe01,WuYe02,WuYe03} and references therein. One exception is  the case where all mappings $F, G, H$ are affine, ${\cal K}$ is polyhedral and $\Lambda$ is the union of finitely many convex polyhedral sets. In this case,  the local error bound property holds automatically following from Robinson's result on polyhedral multifunctions \cite{Robinson81}.
This result, however, depends crucially on the functions  $F, G, H$ being affine and the sets ${\cal K}, \Lambda$  being polyhedral and
the union of finitely many convex polyhedral sets, respectively. The second-order cone, however, is not polyhedral when the dimension is larger than two, and so even when all the mappings $F, G, H$ are affine and the set $\Lambda$ is the  union of  finitely many convex polyhedral sets, the local error bound property may not hold without further assumptions if one of the second-order cones ${\cal K}_i$ has dimension $m_i\geq 3$.
Another easy to verify case is when the gradient vectors $\nabla G_i(z^*) (i=1,\dots, m), \nabla H_i(z^*) (i=1,\dots, m), \nabla F_i(z^*)(i=1,\dots, l)$ are linearly independent. But this condition is very  strong.

The main goal of this paper is to provide a verifiable sufficient condition for  the local error bound property for the second-order cone complementarity system (\ref{system-1})-(\ref{system-2}).  Our condition involves only the first-order and/or the second-order  derivatives of the mappings $F,G,H$   at the point of interest, and is therefore efficiently checkable. The basis of our approach is the sufficient conditions for metric subregularity recently developed by Gfrerer \cite{Gfr11,Gfr13a,Gfr14}, Gfrerer and Klatte \cite{GfrKl16}, Gfrerer and Ye \cite{GfrJYe}. To use these results, we need to compute the tangent cones and the directional normal cones to the second-order cone complementarity set. These results, however, are of independent interest.

{We summarize our main contributions as follows:
\begin{itemize}
\item We introduce a new concept of inner directional normal cone. A set is said to be  directionally regular if the inner directional normal cone coincides with the directional normal cone. It describes the variational geometry of a set along some direction. The directional regularity implies the geometrical derivability. In particular, we show that a convex set is directionally regular. Some useful calculus rules for the directional normal cone are derived.
    \item We establish exact expressions for the tangent cone and the directional normal cone of the second-order cone complementarity set. Moreover we show that the second-order cone complementarity set, which is nonconvex, is directionally regular, and hence both the tangent cone and the directional normal cone commutes with the Cartesian product of finitely many second-order cone complementarity sets.
\item We give sufficient conditions for the existence of error bounds of the second-order cone complementarity problems. These conditions are verifiable based on the initial problem data.
\end{itemize}}
We organize our paper as follows. Section 2 contains the preliminaries.   In Section 3, we study certain properties of the directional   normal cone introduced by  Gfrerer \cite{Gfr13a} and  in Section 4 we derive  sufficient conditions for the  error bound property of a general system by using directional  normal cones.
Section 5 is devoted to the formula and the property of the tangent cone to the second-order complementarity set. In Section 6, we derive the exact expressions for the directional normal cone for the second-order cone complementarity set. Finally in Section 7 we present  sufficient conditions for error bounds of the second-order cone complementarity system.

 The following  notation will be used throughout the paper. We denote by $I$ and $O$ the identity and zero matrix of appropriate dimensions
 respectively. For a matrix $A$, we denote by $A^T$ its transpose.
 The inner product of two vectors $x,y$ is denoted by $x^Ty$ or $\langle x,y \rangle$. For any $z\in \mathbb{R}^n$, we denote by $\|z\|$ the Euclidean norm.
For any nonzero vector $z\in\mathbb{R}^n$, the notation $\bar{z}$ stands for the normalized vector $\frac{z}{\|z\|}$.
 For a function $g: \mathbb{R}^n\rightarrow \mathbb{R}$, we denote $g_+(z):=\max\{0,g(z)\}$, and if it is vector-valued then the maximum is taken componentwise.
{For $z=(z_1,z_2) \in \mathbb{R}\times \mathbb{R}^{m-1}$, we write its reflection about the $z_1$ axis as  $\hat{z}:=(z_1,-z_2)$.}  Denote by $\mathbb{R}z$ the set $\{tz|\ t\in \mathbb{R}\}$. $\mathbb{R}_+z$ and $\mathbb{R}_{++}z$ where $\mathbb{R}_+:=[0,\infty)$ and $\mathbb{R}_{++}:=(0,\infty)$ are similarly defined. For a set ${ C}$,  denote by int$C$, cl$C$, bd$C$, co$C$, $C^c$ its interior, closure, boundary, convex hull, and complement, respectively. The polar cone of a set $C$ is $C^\circ:=\{z|z^Tv\leq 0, \forall v\in C\}$ and $v^\circ$ is the polar cone of a vector $v$. We denote by $d_C(z)$ or $d(z,C)$ the distance from $z$ to $C$. Given a point $z\in \mathbb{R}^n$ and $\varepsilon > 0$, $B_\varepsilon(z)$ denotes an open ball centered at $z$ with radius $\varepsilon$ while $B$ and $\bar{B}$ denote the open and the closed unit ball center at the origin of an appropriate dimension, respectively. For a differentiable  mapping $H:\mathbb{R}^n\to \mathbb{R}^m$ and a vector $z\in \mathbb{R}^n$, we denote by $\nabla H(z)$ the Jacobian matrix of $H$ at $z$.
By $o(\cdot)$, we mean that $o(\alpha)/\alpha \rightarrow 0$ as $\alpha \rightarrow 0$.
 For a set-valued mapping $\Phi: \mathbb{R}^n \rightrightarrows \mathbb{R}^m$,
%(assigning to each $z\in \mathbb{R}^n$ a set $\Phi(z) \subseteq\mathbb{R}^m$ which may be empty),
 the graph and domain of $\Phi$ are denoted by ${\rm gph} \Phi$ and ${\rm dom} \Phi$, respectively, i.e., ${\rm gph}\Phi:=\{(z,v)\in\mathbb{R}^n\times \mathbb{R}^m\,|\ v\in \Phi(z)\}$
%(assigning to each $z\in \mathbb{R}^n$ a set $\Phi(z) \subseteq\mathbb{R}^m$ which may be empty),
and ${\rm dom}\Phi:=\{z\in\mathbb{R}^n\,| \Phi(z)\not =\emptyset\}$. Finally  for any mapping $\varphi:\mathbb{R}^n\rightarrow \mathbb{R}^m$, we denote the active index set at $z^*\in \mathbb{R}^n$ by $I_\varphi(z^*):=\{i\in \{1,\dots,m\}|\varphi_i(z^*)=0\}$. For simplification of notation, we may write $I_\varphi(z^*)$ as $I_\varphi$, provided that there is no confusion in the context.

\section{Preliminaries}
 In this section, we gather some preliminaries on variational analysis and second-order cone which will be used in paper. Detailed discussions on these subjects can be found in \cite{AG03,c,clsw,LVBL98,m1,rw} and the papers we refer to.

\subsection{Background in variational analysis}
Let $\Phi: \mathbb{R}^n \rightrightarrows \mathbb{R}^m$ be a set-valued mapping. We denote by $\limsup_{z'\rightarrow z}\Phi(z')$ and $\liminf_{
z'\rightarrow z}\Phi(z')$ the
 Painlev\'{e}-Kuratowski upper and lower limit, i.e.,
\begin{eqnarray*}
\limsup_{z'\rightarrow z}\Phi(z)&:=&\left  \{v\in \mathbb{R}^m\Big| \exists z_k \rightarrow z, v_k\rightarrow v \mbox{ with } v_k\in \Phi(z_k) \ \ \forall k\right\},\\
\liminf_{z'\rightarrow z}\Phi(z)&:=&\left  \{v\in \mathbb{R}^m\Big| \forall  z_k \rightarrow z, \exists v_k\rightarrow v \mbox{ with } v_k\in \Phi(z_k) \ \  \forall k\right\},
\end{eqnarray*}
respectively.

  Let $C\subseteq \mathbb{R}^n$ and
   $z\in C$. The tangent cone of  $C$
at $ z$ is a closed cone defined by
\begin{equation*}\label{normalcone}
T_C(z)
:=\limsup_{t\downarrow 0}\frac{C-z}{t}
=\left\{u\in\mathbb R^n\Big|\;\exists\,t_k\downarrow
0,\;u_k\to u\;\mbox{ with }\;z+t_k u_k\in C ~\forall  k \right \}.
\end{equation*} The inner  tangent/derivable cone of $C$
at $z$ is defined by
 \[
 T^{i}_C(z):=\liminf\limits_{t\downarrow 0}\frac{C-z}{t}=\left \{u\in\mathbb R^n\Big|\;\forall t_k\downarrow
0,\exists u_k\to u\;\mbox{ such that }\;z+t_k u_k\in C ~\forall  k \right \}.
 \]
 The regular/Fr\'{e}chet normal cone of $C$ at ${z}$ is defined by
  \begin{eqnarray*}
% N^\pi_C ( x^*)&:=& \{v\in \mathbb{R}^n |\  \exists\, M>0\, \mbox{ such that } \langle v, x-x^*\rangle \leq M\|x-x^*\|^2 \ \ %\forall \, x\in C\} \label{eq:def-proximal-norm}\\
 \widehat{{N}}_C({z}) := \Big \{v\in \mathbb{R}^n\,|\
 \langle v, z'-{z}\rangle \leq o(\|z'-z\|) \ \forall  z'\in C \Big\}.\label{normals}
 \end{eqnarray*}
 The limiting/Mordukhovich normal cone is defined by
 % the outer limit of either the proximal normal cone or the regular normal cone, i.e.,
\begin{equation*}\label{eq:def-limit-norm}
N_C (z):=\limsup\limits_{z' \stackrel{C}{\to} z}\widehat{N}_{C}(z')= \Big \{\lim_{i\rightarrow \infty} v_i|\  v_i \in \widehat{N}_C(z_i),\ \   z_i \stackrel{C}{\to} z \Big\}.
\end{equation*}
\begin{defn} We say that a set $C$ is geometrically derivable at a point $z\in C$ if the tangent cone of $C$ coincides with the inner tangent cone of $C$ at $z$, i.e., $T_C(z)=T_C^i(z)$.
\end{defn}

Let $\Phi: \mathbb{R}^n \rightrightarrows \mathbb{R}^m$ be a set-valued mapping and $(x,y) \in {\rm gph } \Phi$. The regular coderivative and the limiting/Mordukhovich coderivative of $\Phi$ at $(x,y)$ are the set-valued mappings defined by
\begin{eqnarray*}
\widehat{D}^*\Phi(x,y)(v)&:=&\left \{ u\in \mathbb{R}^n| (u,-v) \in \widehat{N}_{{\rm gph} \Phi}(x,y)\right \},
\\
{D}^*\Phi(x,y)(v)&:=&\Big \{ u\in \mathbb{R}^n| (u,-v) \in {N}_{{\rm gph} \Phi}(x,y) \Big \},
\end{eqnarray*}
respectively. We omit $y$ in the coderivative notations if the set-valued map $\Phi$  is single-valued at ${x}$.
%, i.e.,  $\Phi(x^*)=\{\bar y\}$.

\medskip

 For a single-valued mapping $\Phi:\mathbb{R}^n \rightarrow \mathbb{R}^m$, the B(ouligand)-subdifferential $\partial_B \Phi$ is defined as
 $$\partial_B \Phi( z)=\left\{\lim_{k\rightarrow \infty} {\nabla} \Phi(z_k)|\ z_k\rightarrow z, \Phi \mbox{ is differentiable at } z_k\right\}.$$
% It is known that ${\rm co}\partial_B \Phi(z)=\partial^c \Phi(z)$, the %Clarke generalized Jacobian of $\Phi$ at $z$. Moreover
If $\Phi$ is a continuously differentiable single-valued map, then
 \[ \widehat{D}^*\Phi(z)={D}^*\Phi(z)=\{{\nabla }\Phi(z)^T\}.
 \]

%Let ${\cal K}$ be the m-dimensional second-order cone. The topological interior and the boundary of ${\cal K}$ are
%\begin{eqnarray*}
%{\rm int} {\cal K}= \{(x_1,x_2) \in \mathbb{R} \times \mathbb{R}^{m-1}|x_1>\|x_2\|\} \ \ {\rm and}\ \  {\rm bd} {\cal K}=\{(x_1,x_2) \in \mathbb{R} \times \mathbb{R}^{m-1}|x_1=\|x_2\|\},
%\end{eqnarray*}
%respectively. For any given vector $x:=(x_1,x_2)\in \mathbb{R}\times \mathbb{R}^{m-1}$, it can be decomposed as
%$$x=\lambda_1(x)c_1(x)+\lambda_2(x)c_2(x),$$
%where $\lambda_i(x)$ and $c_i(x)$ for $i=1,2$ are the spectral values and the associated spectral vectors of $x$ given by
%\[
% \lambda_i(x)=x_1+(-1)^i\|x_2\|  \quad  {\rm and}\quad
% c_i(x)=\left\{\begin{array}{ll}\frac{1}{2}(1, (-1)^i\bar{x}_2) & {\rm if} \ \ x_2\neq 0\\
%\frac{1}{2}(1,w) & {\rm if}\ \ x_2=0\end{array}\right.
%\]
%with $\bar{x}_2:={x_2}/\|x_2\|$ and $w$ being any vector in $\mathbb{R}^{m-1}$ satisfying $\|w\|=1$.
%For $x\in \mathbb{R}^m$, let $\Pi_{\cal K}(x)$ be the metric projection  of $x$ onto ${\cal K}$. Then by \cite{FLT},  it can be calculated by
%\begin{equation}
%\Pi_{\cal K}(x)=(\lambda_1(x))_+c_1(x)+(\lambda_2(x))_+c_2(x).\label{proj}
%\end{equation}

\subsection{Background in variational analysis associated with the second-order cone}
 Let ${\cal K}$ be the $m$-dimensional second-order cone. The topological interior and the boundary of ${\cal K}$ are
\begin{eqnarray*}
{\rm int} {\cal K}= \{(x_1,x_2) \in \mathbb{R} \times \mathbb{R}^{m-1}|\, x_1>\|x_2\|\}, \ \ \ {\rm bd} {\cal K}=\{(x_1,x_2) \in \mathbb{R} \times \mathbb{R}^{m-1}|\, x_1=\|x_2\|\},
\end{eqnarray*}
respectively.

\begin{prop}(see e.g. \cite[Proposition 2.2]{YZ15})\label{Lem3.1}
For any $x,y \in bd {\cal K}\backslash \{0\}$, the following equivalence holds:
$$x^Ty=0 \Longleftrightarrow y=k\hat{x} \mbox{ with } k=y_1/x_1>0  \Longleftrightarrow y=k\hat{x} \mbox{ with } k\in \mathbb{R}_{++}.$$
\end{prop}
%{\bf Proof.} Suppose that $x,y \in bd {\cal K}\backslash \{0\}$ and $x^Ty=0$. Then
% \begin{equation}\label{boundary points-1new}
% x_1=\|x_2\|>0,\ y_1=\|y_2\|>0,\  x^Ty=x_1y_1+x_2^Ty_2=0,
% \end{equation}
% which implies that $-x_2^Ty_2=x_1y_1=\|x_2\|\|y_2\|$.  Hence there exists a positive constant $k$ such that $y_2=-kx_2$.  It follows from (\ref{boundary points-1new}) that $k=y_1/x_1$ and hence  $y=k\hat{x}$. The rest of the proof follows from \cite[Lemma 2.3]{LZL14}.

For any given nonzero  vector $z:=(z_1,z_2)\in \mathbb{R}\times \mathbb{R}^{m-1}$,  we denote by
%$z$ can be  decomposed as
%$$z=\lambda_1(z)c_1(z)+\lambda_2(z)c_2(z),$$
%where $\lambda_1(z),\lambda_2(z)$ and $c_1(z),c_2(z)$ are the spectral values and the associated spectral vectors of $z$ given by
%\begin{eqnarray}
%&& \lambda_1(z)=z_1-\|z_2\|,  \quad  \lambda_2(z)=z_1+\|z_2\|,\label{spva}\\
%%&&
$$c_1(z)=\frac{1}{2}(1, -\bar{z}_2), \quad c_2(z)=\frac{1}{2}(1, \bar{z}_2)
$$ the spectral vectors of $z$, where $\bar{z}_2$ is any vector $w\in \mathbb{R}^{m-1}$ with $\|w\|=1$ if $z_2=0$.
%\label{spve}\end{eqnarray}
%where $\bar{z}_2:=z_2/\|z_2\|$.
For $z\in \mathbb{R}^m$, let $\Pi_{\cal K}(z)$ be the metric projection  of $z$ onto ${\cal K}$ and
% Then by \cite{FLT},  it can be calculated by
%\begin{equation}
%\Pi_{\cal K}(z)=(\lambda_1(z))_+c_1(z)+(\lambda_2(z))_+c_2(z).%\label{proj}
%\end{equation}
$\Pi'_{\cal K}(z;h)$  the directional derivative of $\Pi_{\cal K}$ at $z$ in direction $h$. The following proposition summarizes its formula (see \cite[Lemma 2]{OS08}).

\begin{prop} \label{directd} Let ${\cal K}$ be the $m$-dimensional second-order cone. The mapping $\Pi_{\cal K}(\cdot)$ is directionally differentiable at any $z \in \mathbb{R}^m$ and for any $h\in \mathbb{R}^m$,
\begin{itemize}
\item[{\rm (i)}] if  $z\in int {\cal K}$ or $z\in -int {\cal K}$ or $z\in  (-{\cal K}\cup {\cal K})^c$, then
$\Pi'_{\cal K}(z;h)=\nabla \Pi_{\cal K}(z) h;$
\item[{\rm (ii)}] if $z\in bd {\cal K} \setminus\{0\}$, then $\Pi'_{\cal K}(z;h)=h-2(c_1(z)^Th)_- c_1(z);$
\item[{\rm (iii)}] if $z\in -bd {\cal K} \setminus\{0\}$, then $\Pi'_{\cal K}(z;h)=2(c_2(z)^T h )_+ c_2(z);$
\item[{\rm (iv)}] if $z=0$, then $\Pi'_{\cal K}(z;h)=\Pi_{\cal K}(h)$.
\end{itemize}
\end{prop}
 The following proposition summarizes the regular and the limiting coderivatives of the metric projection operator (see \cite[Lemma 1 and Theorems 1 and 2]{OS08}).

\begin{prop} \label{projection} Let ${\cal K}$ be the $m$-dimensional second-order cone.
\begin{itemize}
\item[{\rm (i)}] If $z\in int {\cal K}$, then $\Pi_{\cal K}$ is differentiable and $\nabla \Pi_{\cal K}(z)=I$.
\item[{\rm (ii)}] If $z\in -int {\cal K}$, then  $\Pi_{\cal K}$ is differentiable and  $\nabla \Pi_{\cal K}(z)=\{O\}$.
\item[{\rm (iii)}] If $z\in  (-{\cal K}\cup {\cal K})^c$, then  $\Pi_{\cal K}$ is differentiable and
 \[
\nabla \Pi_{\cal K} (z)=\frac{1}{2} (1+\frac{z_1}{\|z_2\|})I+\frac{1}{2}\begin{bmatrix}
- \frac{z_1}{\|z_2\|} & \bar{z}_2^T\\
 \bar{z}_2 &  -\frac{z_1}{\|z_2\|}  \bar{z}_2 \bar{z}_2^T
 \end{bmatrix}.
 \]
\item[\rm (iv)] If $z\in bd {\cal K} \setminus\{0\}$, then
\begin{eqnarray*}
\widehat{D}^*\Pi_{\cal K}(z)(u^*)&=&\{x^*|u^*-x^*\in \mathbb{R}_+c_1(z), \langle x^*, c_1(z)\rangle \geq 0\},\\
{D}^*\Pi_{\cal K}(z)(u^*)&=&{\partial}_B \Pi_{\cal K}(z)u^* \cup \{x^*|u^*-x^*\in \mathbb{R}_+c_1(z),  \langle x^*, c_1(z)\rangle \geq 0\},
\end{eqnarray*}
%where $c_1(z):=\frac{1}{2} (1,-\bar{z}_2)$ a
and
$${\partial}_B \Pi_{\cal K}(z)=\left\{ I,I+\frac{1}{2}\begin{bmatrix}
- 1 &  \bar{z}_2^T \\
 \bar{z}_2 &  - \bar{z}_2 \bar{z}_2^T
 \end{bmatrix} \right\}.$$
\item[\rm (v)] If $z\in -bd {\cal K} \setminus\{0\}$, then
\begin{eqnarray*}
\widehat{D}^*\Pi_{\cal K}(z)(u^*)&=&\{x^*|x^*\in \mathbb{R}_+c_2(z), \langle u^*-x^*, c_2(z)\rangle \geq 0\}, \\
{D}^*\Pi_{\cal K}(z)(u^*)&=&{\partial}_B \Pi_{\cal K}(z)u^* \cup \{x^*|x^*\in \mathbb{R}_+c_2(z),  \langle u^*-x^*, c_2(z)\rangle \geq 0\},
\end{eqnarray*}
%where   $c_2(z)=\frac{1}{2} (1,\bar{z}_2)$
 and
$${\partial}_B \Pi_{\cal K}(z)=\left\{O,\frac{1}{2}\begin{bmatrix}
1 &  \bar{z}_2^T\\
 \bar{z}_2 &  \bar{z}_2 \bar{z}_2^T
 \end{bmatrix} \right\}.$$
\item[\rm (vi)] If $z=0$, then
\begin{eqnarray*}
\widehat{D}^*\Pi_{\cal K}(z)(u^*)&=&\{x^*|\ x^* \in {\cal K},\ u^*-x^* \in {\cal K}\}.\\
{D}^*\Pi_{\cal K}(z)(u^*)&=&{\partial}_B \Pi_{\cal K}(0)u^* \cup \{ x^*|\, x^* \in {\cal K},\ u^*-x^*\in {\cal K}\}\\
& &\cup  \bigcup_{\xi \in C} \{x^*|\, u^*-x^*\in \mathbb{R}_+\xi,\  \langle x^*, \xi\rangle \geq 0\}\\
& & \cup \bigcup_{\eta \in C} \{x^*|\, x^*\in \mathbb{R}_+\eta, \  \langle u^*-x^*, \eta \rangle \geq 0\},
\end{eqnarray*}
where $C:=\{\frac{1}{2}(1,w)|\ w\in \mathbb{R}^{m-1},\ \|w\|=1\}$ and
$$\partial_B \Pi_{\cal K}(0)=\{O,I\}\cup \left\{\left.\frac{1}{2}\begin{bmatrix}
 1 & w^T \\
 w & 2\alpha I+(1-2\alpha)ww^T
 \end{bmatrix}\right|\, w\in \mathbb{R}^{m-1}, \ \ \|w\|=1, \ \ \alpha\in [0,1] \right\}.
 $$
\end{itemize}
\end{prop}

\begin{prop}\cite[Proposition 2.1]{yzhou} \label{Prop2.2} Let $(x,y)\in \Omega:=\{(x,y)|x\in {\cal K}, y\in {\cal K}, x^Ty=0\}$. Then
\begin{eqnarray*}
 \widehat{N}_\Omega (x,y)&=&\bigg\{(u,v)|\ -v\in \widehat{D}^*\Pi_{\cal K}(x-y)(-u-v)\bigg\},\\
 {N}_\Omega (x,y)&=&\bigg\{(u,v)|\ -v\in {D}^*\Pi_{\cal K}(x-y)(-u-v)\bigg\}.
 \end{eqnarray*}
 \end{prop}

The exact formula of the regular normal cone and limiting normal cone of $\Omega$ have been
 established in \cite{yzhou}.

 \begin{prop}\cite[Theorem 3.1]{yzhou}\label{formula-regular normal conenew}
 Let $(x,y)$ be  in the $m$-dimensional second-order cone complementarity set $\Omega$. Then
 \begin{eqnarray*}
 \widehat{N}_\Omega(x,y)=&   \left  \{ \begin{array}{ll}
\{(u,v)| u\in \mathbb{R}^m, \ v=0 \} & \ {\rm if} \ x=0,\ y\in {\rm int}{\cal K}; \\
\{(u,v)| u=0, v\in \mathbb{R}^m  \}  & \ {\rm if}\ x\in {\rm int}{\cal K}, \ y=0 ;\\
\{(u,v)|  u\perp x , \ v\perp y,  \ x_1\hat{u}+y_1v\in \mathbb{R} x \}   & \ {\rm if}\ x,y \in {\rm bd}{\cal K}\backslash\{0\} , x^Ty=0; \\
\{(u,v)| u\in \hat{y}^\circ,\ v\in \mathbb{R}_-\hat{y}  \} &  \ {\rm if}\ x=0, \ y\in {\rm bd}{\cal K}\backslash \{0\} ;\\
 \{(u,v)|  u\in \mathbb{R}_-\hat{x}, v\in \hat{x}^\circ   \}  & \ {\rm if}\ x\in {\rm bd}{\cal K}\backslash \{0\},\ y=0 ;\\
\{(u,v)|   u\in -{\cal K},\ v\in -{\cal K}  \}  & \ {\rm if}\ x=0, \ y=0.
\end{array}\right.
 \end{eqnarray*}
 \end{prop}

 \begin{prop}\cite[Theorem 3.3]{yzhou} \label{formula-normal cone}
   Let $(x,y)$ be in the $m$-dimensional second-order cone complementarity set $\Omega$. Then
 \begin{eqnarray*}
 N_\Omega(x,y)= \widehat{N}_\Omega (x,y)=\left\{ \begin{array}{ll}
 \{(u,v)| u\in \mathbb{R}^m, \ v=0 \}  &\ {\rm if} \ x=0,\ y\in {\rm int}{\cal K}; \\
\{(u,v)| u=0, \ v\in \mathbb{R}^m \} & \ {\rm if}\ x\in {\rm int}{\cal K},\ y=0 ;\\
 \{(u,v)| u\perp x,\ v\perp y,\  x_1\hat{u}+y_1v\in \mathbb{R} x\}  &\ {\rm if} \  x, y\in {\rm bd}{\cal K}\backslash\{0\}.
 \end{array} \right.
 \end{eqnarray*}
 For $x=0, y\in {\rm bd}{\cal K}\backslash \{0\}$,
 $$ N_\Omega(x,y)=\{(u,v)|u\in \mathbb{R}^m,\ v=0
  \ \ {\rm or} \ \
 u\perp \hat{y},\ v\in \mathbb{R} \hat{y}
 \ \ {\rm or} \ \ \langle u, \hat{y} \rangle \leq 0, \ v\in \mathbb{R}_- \hat{y}\};$$
 for $x\in {\rm bd}{\cal K}\backslash \{0\}, y=0$,
 $$ N_\Omega(x,y)=\{(u,v)|  u=0, v\in \mathbb{R}^m \ \ {\rm or}\ \ u\in \mathbb{R} \hat{x},\
 v\perp \hat{x}\ \ {\rm or} \ \ u\in \mathbb{R}_-\hat{x},
  \langle v, \hat{x} \rangle \leq 0 \};$$
  for $x=y=0$,
 \begin{eqnarray*}
N_\Omega(x,y)&=&  \{(u,v)|\  u\in -{\cal K}, v\in -{\cal K} \ \mbox{\rm or }  u\in \mathbb{R}^m, v=0 \
 \mbox{\rm or } u=0, v\in \mathbb{R}^m\\
  &&  \  \mbox{\rm or } u\in \mathbb{R}_-\xi,\ v\in \xi^\circ \ \  \mbox{\rm or }\ u\in \xi^\circ, v\in \mathbb{R}_-\xi\\
   && \ \mbox{\rm or } u\perp \xi, \ v\perp \hat{\xi}, \
 \alpha \hat{u}+(1-\alpha)v\in \mathbb{R} \xi, \ \mbox{ for some }
  \alpha\in [0,1],  \xi \in C\}
 \end{eqnarray*}
 where \[
 C:=\{(1,w)|\ w\in \mathbb{R}^{m-1},\ \|w\|=1\}.\label{C}\]
  \end{prop}

 \section{Calculus for directional normal cones}
 Recently
a directional version of the limiting normal cone was  introduced by  Gfrerer \cite{Gfr13a} and used to derive sufficient conditions for metric subregularity, which form the basis for our approach. Since   calculus for the directional normal cone is very important and the existing results are rather rare, in this section, we develop some calculus for the directional normal cone. First, we recall the definition of a directional normal cone.
\begin{defn}Given a set
$C\subseteq\mathbb R^n$, a point $z\in C$
and  a direction $d\in \mathbb{R}^{n}$, the limiting normal cone to $C$ in  direction $d$ at $z$ is defined by
\[
N_{C}(z; d):=\limsup\limits_{{t\downarrow 0}\atop {d'\to d}}\widehat{N}_C(z+td')=\left \{v | \exists t_{k}\downarrow 0, d_{k}\rightarrow d, v_{k}\rightarrow v \ {\rm with }\  v_{k}\in \widehat{N}_{C}(z+ t_{k}d_{k}), \ \forall k \right \}.
\]
\end{defn}
%It follows from the definition that $N_{C}(z; d)=\emptyset$ for any $d\not \in T_C(x)$. The following proposition gives the connection between the directional normal cone and the limiting normal cone in the case when the set $C$ is a union of finitely many closed convex sets.
%\begin{prop}\cite[Lemma 2.1]{Gfr14}\label{Lemma2.1} Let $C\subset \mathbb{R}^n$ be the union of finitely many closed convex sets $C_i, i=1,\dots, m$, $ x \in C, u\in \mathbb{R}^n$. Then
%$$N_C(x;u) \subseteq\{v\in N_C(x)| v^Tu=0\}.$$
%If moreover $C$ is convex and $u\in T_C(x)$, then $N_C(x;u) = N_C(x)$.
%\end{prop}

We define the concept of the inner directional normal cone as follows.
 \begin{defn}
Given a set
$C\subseteq\mathbb R^n$, a point $z\in C$
and  a direction $d\in \mathbb{R}^{n}$, the inner limiting normal cone to $C$ in  direction $d$ at $z$ is defined by
\[
N^i_{C}(z;d):=\left \{v|\ \forall t_k\to 0, \exists d_k\to d, v_k\to v \ {\rm with}\ v_k\in \widehat{N}_{C}(z+t_kd_k), \ \forall k \right \}.
\]
\end{defn}
The following results follow from definition immediately.
\begin{prop}\label{aa-co-1} For any set $C$ and any $z\in C$,
\[
{\rm dom}N_C(z;\cdot)=T_C(z), \ \ {\rm dom}N^i_C(z;\cdot)=T^i_C(z).
\]
\end{prop}
%\beginproof
%It is clear that $N^{i}_{A}(x;d)=\emptyset$ if $d\notin T^{i}_A(x)$. If $d\in T^{i}_A(x)$, then $0\in \widehat{N}_{A}(x+t_nd_n)$, so $0\in N^i_{A}(x;d)$. The argument can be applicable to $N_A(x;d)$.
%\endproof
{It is easy to see that} $N^i_{C}(z;d)\subseteq N_C(z;d)\subseteq N_C(z)$ for any $d$ and
  $N_{C}(z; 0)=N_{C}(z)$.
\begin{defn}
{Given a  subset $C$ in $\mathbb{R}^n$ and $z\in C, d\in \mathbb{R}^n$, we say that the set $C$ is regular at $z$ in direction $d$ if
\begin{equation}\label{directional-2}
N_{C}(z;d)=N^i_{C}(z;d).
\end{equation}
If the above formula holds for all $d$, we say that $C$
is directionally regular at $z$. If $C$ is directionally regular at any point ${z\in C}$, then we say that $C$ is directionally regular.}
\end{defn}
It is clear that set $C$ is regular in direction $d$ for any $d\notin T_C(z)$, since both sides of (\ref{directional-2}) are empty in this case.
It follows from Proposition \ref{aa-co-1}  that the directional regularity implies the geometric derivability.
\begin{cor}\label{geomd}
If the set $C$ is directionally regular at $z\in C$, then $C$ is geometrically derivable at $z\in C$.
\end{cor}

An important property of the limiting normal cone is that it commutes with the Cartesian product (see e.g.
\cite[Proposition 1.2]{m1}): for any sets $A_1,\dots, A_I$,
\begin{eqnarray*}
N_{A_1\times\cdots\times A_I}(z_1,\dots,z_I)=N_{A_1}(z_1)\times \cdots \times N_{A_I}(z_I).
\end{eqnarray*} It is easy to verify that this property holds for the inner directional normal cones, i.e.,
\begin{equation}\label{product-inner normal cone}
N^i_{A_1\times\cdots\times A_I}((z_1,\dots,z_I);(d_1,\dots, d_I))=N^i_{A_1}(z_1;d_1)\times \cdots \times N^i_{A_I}(z_I;d_I).
\end{equation}
For directional normal cones, this kind of property does not come free. Fortunately, it holds under the directional regularity.

\begin{prop}\label{product2}
The inclusion
\[
N_{A\times B}((x,y);(d,w))\subseteq N_{A}(x;d)\times N_{B}(y;w)
\] holds for any given sets $A\subseteq \mathbb{R}^n ,B\subseteq
\mathbb{R}^m$,  any point $(x,y)\in A\times B$, and any direction $(d,w)\in \mathbb{R}^n\times \mathbb{R}^m$. Moreover if either $A$  is
 regular at $x$ in direction $d$ or $B$ is regular at $y$ in direction $w$, then
\[
N_{A\times B}((x,y);(d,w))=N_{A}(x;d)\times N_{B}(y;w).
\]
If $A$ and $B$ are regular in directions $d,w$ respectively, then $A\times B$ is  regular in direction $(d,w)$.
\end{prop}

\beginproof Note that
\begin{eqnarray*}
N_{A\times B}((x,y);(d,w))&=& \limsup\limits_{{t\downarrow 0}\atop {(d',w')\to (d,w)}}\widehat{N}_{A\times B}\big((x,y)+t(d',w')\big)\\
&=&  \limsup\limits_{{t\downarrow 0}\atop {(d',w')\to (d,w)}}\widehat{N}_{A}(x+td')\times \widehat{N}_{B}(y+tw')\\
&\subseteq&  \limsup\limits_{{t\downarrow 0}\atop {d'\to d}}\widehat{N}_{A}(x+td')\times \limsup\limits_{{t\downarrow 0}\atop {w'\to w}} \widehat{N}_{B}(y+tw')\\
&=& N_{A}(x;d)\times N_{B}(y;w).
\end{eqnarray*}
Conversely, take $(p,q)\in N_{A}(x;d)\times N_{B}(y;w)$. Without loss of generality, assume that $A$ is regular at $x$ in direction $d$. Since $q\in N_{B}(y;w)$, there exists $t_n\downarrow 0$ and $w_n\to w$ and $q_n\to q$ such that $q_n\in \widehat{N}_{B}(y+t_nw_n)$. Since $p\in N_{A}(x;d)$ and $A$ is regular at $x$ in direction $d$,  for the above $t_n$ there exist $d_n\to d$ and $p_n\to p$ such that $p_n\in \widehat{N}_{A}(x+t_nd_n)$. So
\[
(p_n,q_n)\in \widehat{N}_{A}(x+t_nd_n)\times \widehat{N}_{B}(y+t_nw_n)=\widehat{N}_{A\times B}((x,y)+t_n(d_n,w_n)).
\]
By definition, this means that $(p,q)\in N_{A\times B}((x,y);(d,w))$.

Now suppose that $A$ and $B$ are regular in direction $d,w$ respectively. Then
\begin{eqnarray*}
N_{A\times B}((x,y);(d,w)) &\subseteq &N_{A}(x;d)\times N_{B}(y;w)\\
& =& N^i_{A}(x;d)\times N^i_{B}(y;w)\\
&=& N^i_{A\times B}((x,y);(d,w)).
\end{eqnarray*}
Therefore $A\times B$ is regular in direction $(d,w)$.
\endproof

 \begin{prop}\label{TanNor}
Let $(z_1,\dots,z_I)\in A_1\times \cdots \times A_I$ and $(d_1,\dots,d_I)$ be given. Then
\begin{eqnarray}
&& T_{A_1\times\cdots\times A_I}(z_1,\dots,z_I)\subseteq T_{A_1}(z_1)\times \cdots \times T_{A_I}(z_I),\label{tangent-product}\\
&&
N_{A_1\times\cdots\times A_I}\big((z_1,\dots,z_I);(d_1,\dots,d_I)\big) \subseteq N_{A_1}(z_1;d_1)\times \cdots \times N_{A_I}(z_I;d_I),\label{directional-product}
\end{eqnarray}
and equality holds if all except at most one of $A_i$ for $i=1,\dots,I$ are directionally regular at $z_i$.
 \end{prop}
 \beginproof
By Corollary \ref{geomd}, the directional regularity implies the geometric derivability. Then the tangent set formula follows from applying \cite[Proposition 1]{GfrJYe}.
 The directional normal cone formula follows from Proposition \ref{product2}.
 \endproof

In the rest of this section we will study some calculus rule of the directional normal cone, and in the mean time examine the directional regularity.

 \begin{prop}\label{lem-directional}
 If $C\subseteq\mathbb R^n$ is a closed cone, then
 \[N^i_C(0;d)=N_{C}(0;d)=N_{C}(d),  \ \forall d\in \mathbb{R}^n.\]
\end{prop}
\beginproof
 Since $C$ is a cone, we have $T_{C}(0)=C$. If $d\not \in C$, then \[{N_C^i(0;d)=}N_{C}(0;d)=N_{C}(d)=\emptyset.\] If $d\in C$, then
\[N_{C}(0;d)=\limsup\limits_{t\downarrow 0\atop d'\to d}\widehat{N}_{C}(td')=
\limsup\limits_{d'\to d}\widehat{N}_{C}(d')=N_{C}(d).\]
%Note that
%\begin{equation}\label{aa-9-4}
%\widehat{N}_{C}(x)=\widehat{N}_{C}(\lambda x), \ \ \forall \lambda>0.
% \end{equation}
% In fact, if $v\in \widehat{N}_C(x)$, then
%\begin{equation}\label{aa-9-3}
%\limsup\limits_{z\in C, z\to x}\frac{\langle v, z-x \rangle }{\|z-x\|}\leq 0
%\end{equation}
%Then
%\[
%\limsup\limits_{z\in C, z\to \lambda x}\frac{\langle v, z-\lambda x \rangle }{\|z-\lambda x\|}=
%\limsup\limits_{\frac{z}{\lambda}\in C, \frac{z}{\lambda} \to x}\frac{\langle v, \frac{z}{\lambda}-x \rangle }{\|\frac{z}{\lambda}- x\|}
%\leq 0
%\]
%where the last step is due to (\ref{aa-9-3}). So $v\in \widehat{N}_C(\lambda x)$. Then $\widehat{N}_C(x)\subseteq\widehat{N}_C(\lambda x)$. The converse inclusion can be shown similarly.
Now we show that $N^i_C(0;d)=N_C(0;d)$. It suffices to show $N_C(0;d)\subseteq N^i_C(0;d)$. Take $v\in N_C(0;d)$. Then there exists $\eta_n\downarrow 0$ and $d_n\to d$ and
$v_n\to v$ with $v_n\in \widehat{N}_C(\eta_n d_n)$. For any $t_n\downarrow 0$, take $d_n$ and $v_n$ above, then $d_n\to d$ and $v_n\to v$ with $v_n\in \widehat{N}_C(\eta_n d_n)=\widehat{N}_C( d_n)=\widehat{N}_C(t_n d_n)$. Hence
$v\in N^i_C(0;d)$.
\endproof

We next show that any convex set {is regular along any direction}.

\begin{prop}\label{convex}
Any closed convex set $A$ is directionally regular.
\end{prop}

\beginproof
Since $N^i_{A}(z;d)\subseteq N_A(z;d)$ for any $z$ and $d$, it suffices to prove $N_{A}(z;d)\subseteq N^{i}_{A}(z;d)$ for any $d\in T_A(z)$. Take $w\in N_{A}(z;d)$ with $d\in T_A(z)$. Then there exists $\eta_k\downarrow 0$, $d_k\to d$ and $w_k\to w$ with $w_k\in \widehat{N}_{A}(z+\eta_kd_k)$. Since $A$ is convex, it follows that
\begin{equation}\label{a9-1}
\langle w_k, z'-z-\eta_kd_k  \rangle\leq 0, \ \  \ \ \forall z'\in A.
\end{equation}
In particular, taking $z'=z$ in the above, we have
\begin{equation}\label{a9-2}
\langle w_k, d_k\rangle\geq 0.
\end{equation}
%We now show that $w\in N^{i}_A(x;d)$.
Let  $t_n\downarrow 0$. Then since $\eta_k\downarrow 0$,  for each fixed $n$, there exists $k(n)$ satisfying $k(n)\geq n$ and $\eta_{k(n)}<t_n$.
   Hence $k(n)\to \infty$ as $n\to \infty$. For simplicity, denote by $d_n:=d_{k(n)}$ and $w_n:=w_{k(n)}$. Since $\{d_n\}$ and $\{w_n\}$  are subsequences of $\{d_k\}$ and $\{w_k\}$ respectively, we have
$d_n\to d$ and $w_n\to w$. Hence for all $z'\in A$ we have
\begin{eqnarray*}
\langle w_{n}, z'-z-t_nd_n  \rangle &=& \langle w_{k(n)}, z'-z-t_nd_{k(n)}  \rangle \\
&=& \langle w_{k(n)}, z'-z-\eta_{k(n)}d_{k(n)}  \rangle+
\langle w_{k(n)}, \eta_{k(n)}d_{k(n)}-t_nd_{k(n)}\rangle\\
&\leq & \langle w_{k(n)}, \eta_{k(n)}d_{k(n)}-t_nd_{k(n)}\rangle\\
&\leq & 0,
\end{eqnarray*}
where the first inequality comes from (\ref{a9-1}) and the second inequality follows from (\ref{a9-2}).
So $w_n\in N_{A}(z+t_nd_n)=\widehat{N}_A(z+t_nd_n)$. By the definition of the inner limiting normal cone, we have $w\in N^{i}_A(z;d)$. This completes the proof.
\endproof

%\textcolor{red}{Can we prove the following proposition? See \cite[Lemma 1]{GfrKl16}}
%\begin{prop}
%Let $A$ be is a union of finitely many polyhedral sets. Then $A$ is directional normally regular.
%\end{prop}

{
Based on (\ref{product-inner normal cone}), Propositions \ref{TanNor} and \ref{convex}, we can obtain the following results.}

\begin{cor}\label{corollary-product} Let $A_i$ be given for $i=1,\dots,I$.
\begin{enumerate}
\item[\rm (i)] If $A_i$ is regular at $z_i\in A_i$ in direction $d_i$ for $i=1,\dots,I$, then $ A_1\times\dots\times A_I$ is regular at  $(z_1,\dots,z_I)$ in direction $(d_1,\dots,d_I)$. Moreover, (\ref{directional-product}) holds as an equation.
\item[\rm (ii)] If $A_i$ is directionally regular at $z_i\in A_i$ for $i=1,\dots,I$, then  $A_1\times \cdots \times A_I$ is directionally regular at $(z_1,\dots,z_I)$. Moreover (\ref{tangent-product}) and (\ref{directional-product}) holds as an equation for all $d_i$.
\item[\rm (iii)] If $A_i$ is closed and convex for $i=1,\dots,I$, then $A_1\times \cdots \times A_I$ is directionally regular. Moreover (\ref{tangent-product}) and (\ref{directional-product}) holds as an equation for all $z_i\in A_i$ and $d_i$.
\end{enumerate}
\end{cor}

  Corollary \ref{corollary-product} (iii) extends the result given in \cite[Lemma 1]{GfrKl16}, where each $A_i$ is assumed to be a polyhedral convex set.
  In Section 6, we will show that the second-order cone complementarity set, although it is a nonconvex set, is directionally regular.

\section{Sufficient conditions for the error bound property via  directional normal cones}
Consider a general system in the form:
$ P(z) \in D,$
where $P:\mathbb{R}^l \rightarrow \mathbb{R}^s$  and $D\subseteq \mathbb{R}^s $  is closed.
%\begin{defn}
%  	 Let $M:\mathbb{R}^d   \rightrightarrows \mathbb{R}^s$ be a set-valued mapping and let $(z^*,0)\in gph M$. We say that the system $0\in M(z)$ has a local error bound at $z^*$ or $M$ is {\em metrically subregular}  at $(z^*,0)$ if there exist a neighborhood $V$ of $z^*$ and a positive number  $\kappa>0$ such that
%  \begin{equation}\label{EqMetrSubReg}d(z,M^{-1}(0))\leq\kappa d(0,M(z))\ \; \forall z\in V.
%  \end{equation}
%\end{defn}
We say that the system $ P(z) \in D$ has a local error bound at $z $ such that $P(z)\in D$, or the set-valued mapping $M(z):=P(z)-D$ is metrically subregular at $(z,0)\in {\rm gph} M$,   if there exist a neighborhood $V$ of $z$ and a positive number  $\kappa>0$ such that
 \[d_{M^{-1}(0)}(z')\leq\kappa d_D(P(z')),\ \; \forall z'\in V.
  \]
It is easy to see that $M$ is metrically subregular at  $(z,0)$ if and only if its inverse set-valued mapping $M^{-1}$ is  calm at $(0,z)\in {\rm gph} M^{-1}$, i.e.,   there exist a neighborhood $W$ of $0$,  a neighborhood $V$ of $z$ and a positive number  $\kappa>0$ such that
  $$M^{-1}(w)
  \cap V\subseteq
   M^{-1}(0) +
   \kappa
   \|w\| {\bar B}, \quad \forall w\in W. $$
   The metric subregularity is obviously weaker than the metric regularity (or the pseudo Lipschitz continuity) which ensures the
 existence of  a neighborhood $W$ of $0$,  a neighborhood $V$ of $z$ and a positive number  $\kappa>0$ such that
  $$M^{-1}(w)
  \cap V\subseteq
   M^{-1}(w') +
   \kappa
   \|w-w'\| {\bar B}, \quad \forall w,w'\in W. $$

     While the term for the calmness of a set-valued map was first coined in \cite{rw}, it was introduced as the pseudo-upper Lipschitz continuity  in \cite{yy}, taking into the account that it is weaker than both the pseudo-Lipschitz continuity of Aubin \cite{Aub} and the upper Lipschitz continuity of Robinson \cite{Robinson75,Robinson76}. More information and discussion on metric regularity and the related concept can be found in \cite{KK02}.

 Recall that the following well-known criteria for metric regularity of the set-valued mapping {$M$ or the Aubin property of its inverse mapping $M^{-1}(w)=\{ z\in \mathbb{R}^l| P(z){-w}\in D\}$}.

   \begin{thm}(see e.g. \cite{rw}) \label{Thm4.1} Consider the system $P(z)\in D$, where $P$ is smooth and $D$ is closed. Then the set-valued map {$M(z):=P(z)-D$ is metrically regular at $(z,0)$} if and only if the no nonzero abnormal multiplier constraint qualification (NNAMCQ) holds at $z$,  i.e.,
       \begin{equation}\nabla P(z)^T\lambda=0,\;\lambda\in N_D(P(z)) \quad\Longrightarrow\quad\lambda=0.\label{Mordukhovich-1}  \end{equation}
%Conversely the condition NNAMCQ is necessary for metric regularity %provided that the set $D$ is convex.
 \end{thm}
  While following \cite{ye}, the  condition  (\ref{Mordukhovich-1}) is called NNAMCQ, there are other terminologies in the literature; e.g.,   generalized MFCQ (GMFCQ) in \cite{FKO07} and Mordukhovich criterion in \cite{GfrKl16}. This condition is a necessary and sufficient condition for metric regularity and hence may be too strong for metric subregularity.

By using the directional normal cone instead of the limiting normal cone, the following sufficient conditions for metric subregularity have been introduced.
   \begin{thm}(\cite[Corollary 1]{GfrKl16}) \label{ThSuffCondMS}Let  $ P(z)\in D$ with $P$ smooth.
%Let $\zb\in\Omega$ be feasible for \eqref{EqGenOptProbl}.
The set-valued mapping  $M(z):=P(z)-D$ is metrically subregular at $(z,0)$ if the
 first-order sufficient condition for metric subregularity (FOSCMS) holds:  for every $0\not=w$ such that $\nabla P(z)w\in T_D(P(z))$ one has
      \[\nabla P(z)^T\lambda=0,\;\lambda\in N_D(P(z);\nabla P(z)w)\quad\Longrightarrow\quad\lambda=0.\]
\end{thm}

Let us discuss the relation between FOSCMS and NNAMCQ. {FOSCMS can be rewritten equivalently as
\begin{equation}\label{FOSCMS-2}
\nabla P(z)^T\lambda=0, \ \  \lambda\in \bigcup\limits_{w\in \Gamma}N_{D}(P(z);\nabla P(z)w)\Longrightarrow \lambda=0,
\end{equation}
where $\Gamma:=\{w\neq 0
|\nabla P(z)w\in T_{D}(P(z))\}$. Noth that FOSCMS holds automatically if $\Gamma=\emptyset$, i.e.,
\begin{equation}\nabla P(z)w\in T_{D}(P(z)) \Longrightarrow w=0.\label{eqnnew1}
\end{equation}}
%\noindent either
%\begin{equation}\nabla P(z)w\in T_{D}(P(z)) \Longrightarrow w=0\label{eqnnew1}
%\end{equation}
%or \begin{equation}
%\exists 0\neq w \mbox{ s.t. } \nabla P(z)w\in T_{D}(P(z))\label{eqnnew2}
%\end{equation} such that
%  \begin{equation}\label{FOSCMS-2}
%\nabla P(z)^T\lambda=0, \ \  \lambda\in \bigcup\limits_{\{0\neq w
%|\nabla P(z)w\in T_{D}(P(z))\}}N_{D}(P(z);\nabla P(z)w)\Longrightarrow \lambda=0.
%\end{equation}
%where $T^{\rm lin}_{\Gamma}(z):=\{w|\nabla P(z)w\in T_{D}(P(z))\}$.
According to the graphical derivative criterion for strong metric subregularity \cite{DoRo14}, condition (\ref{eqnnew1})  is equivalent to saying that the set-valued map $M(z)=P(z)-D$ is strongly metrically subregular (or equivalently its inverse is isolatedly calm) at $( z,0)$.

\begin{thm}\label{thm4.3}
Let $M(z):=P(z)-D$ and $(z,0)\in {\rm gph}M$. { FOSCMS at $z$} is equivalent to NNAMCQ at $z$ under one of the following assumptions:
\begin{itemize}
\item[\rm (i)]
$\nabla P(z)$ does not have full column rank;
\item[\rm (ii)] $D$ is a closed and convex set and there exists $
 w\neq 0$ such that $ \nabla P(z)w\in T_{D}(P(z)).$
\end{itemize}
\end{thm}

\beginproof
 (i). If $\nabla P(z)$ does  not have full column rank, then there exists $\bar{w}\neq 0$ such that $\nabla P(z)\bar{w}=0$.  So $\Gamma\neq \emptyset$.
 %Hence (\ref{eqnnew1}) does not hold and FOSCMS means that (\ref{eqnnew2}) and (\ref{FOSCMS-2}) hold.
 Since
 \[
 N_{D}(P(z);\nabla P(z)\bar{w})=N_{D}(P(z);0)=N_{D}(P(z)),
 \]
 we have
 \[
 \bigcup\limits_{w\in \Gamma}N_{D}(P(z);\nabla P(z)w)=N_{D}(P(z)).\]
 It follows that FOSCMS and NNAMCQ are equivalent by comparing the conditions (\ref{Mordukhovich-1})
 and (\ref{FOSCMS-2}).
  %Conversely, suppose NNAMCQ holds.  In this case, (\ref{eqnnew1}) %does not hold. So by virtue of  $N_{D}(z; d)\subset N_{D}(z)$, we %have FOSCMS holds.

 (ii). Suppose that $D$ is a closed and convex set and there exists $w\neq 0$ such that  $\nabla P(z)w\in T_{D}(P(z))$.   Then FOSCMS means (\ref{FOSCMS-2}) holds. Since the directional normal cone is in general a subset of the limiting normal cone, it is clear that NNAMCQ implies FOSCMS. Conversely assume that FOSCMS holds.
 %Then  there exists $\bar{u}\neq 0$ such that $\nabla P(\bar{z})%\bar{u}\in T_{D}(P(\bar{z})$.
  Take $\lambda$ satisfying
 $ \nabla P(z)^T\lambda=0$ and $\lambda\in N_{D}(P({z}))$.
% If there does not exist nonzero $u$ satisfying  $\nabla P(\bar{z})u\in T_{D}(P(\bar{z}))$, then $M$ has strong metric subregularity, and hence is metric regular. If there exists $0\neq u$ satisfying $\nabla P(\bar{z})u\in T_{D}(P(\bar{z}))$,
 Note that
 $\langle \lambda, \nabla P({z})w \rangle=\langle \nabla P({z})^T\lambda, w \rangle=0$. Hence $\lambda\in N_{D}(P({z}))\cap (\nabla P({z})w)^\perp$, which means that $\lambda \in N_{D}(P({z});\nabla P({z})w)$ by \cite[Lemma 2.1]{Gfr14}. The FOSCMS at ${z}$ then ensures $\lambda=0$. Hence NNAMCQ holds at ${z}$.
 \endproof

 \begin{remark} The assumptions (i) or (ii) given in Theorem \ref{thm4.3} cannot be omitted. For example, when $\nabla P(z)$ has full column rank and $D$ is nonconvex with $\Gamma\neq \emptyset$,  FOSCMS may be strictly weaker than NNAMCQ; see  Example \ref{ex4.1} below. \end{remark}
\begin{example}\label{ex4.1} Consider the optimization problem:
 \[\begin{array}{lll}
 &\min & z_1+z_2\\
 &\ s.t. & (z_1,z_2)\in \K\\
 && (z_1,z_2)\in \Omega,
 \end{array}\]
 where ${\cal K}:=\{(z_1,z_2)\in \mathbb{R}^2|z_1\geq |z_2|\}$ and
 $\Omega:=\{(z_1,z_2)\in \mathbb{R}^2| z_1\geq 0 , z_2\geq 0, z_1 z_2=0\}$.
Denote by  $P(z)=(z,z)$ and $D=\K\times \Omega$. The optimal solution is $z^*=(0,0)$.  It is clear that \[\{w\neq 0| \nabla P(z^*) w\in T_{D}(P(z^*))\}
 =\{w\neq 0| (w,w)\in D\}=\{(w_1,w_2)| w_1>0, w_2=0\},\]  $\nabla P(z^*)$ has a
 full column rank and $D$ is nonconvex.  By virtue of Proposition \ref{lem-directional}, since $D$ is a cone we have $N_{D}((0,0);\nabla P(z^*)w)=N_{D}(\nabla P(z^*)w)$, and hence the condition
 \[
 \nabla P(z^*)^T\lambda=0, \ \ \lambda\in N_{D}(P(z^*);\nabla P(z^*)w)=N_{D}(\nabla P(z^*)w)=N_{\cal K}(w_1,w_2)\times N_\Omega (w_1,w_2)
 \]
 with $w_1>0, w_2=0$
 takes the form
 \[
 \lambda^{\K}+\lambda^{\Omega}=0, \ \ (\lambda^{\K}, \lambda^{\Omega})\in \{(0,0)\}\times (\{0\}\times \mathbb{R}),
 \]
which implies that $(\lambda^{\K},\lambda^{\Omega})=(0,0)$. Hence FOSCMS holds at $z^*$.

 On the other hand,
  \[
 \nabla P(z^*)^T\lambda=0, \ \ \lambda=(\lambda^{\cal K}, \lambda^\Omega)\in N_{D}(P(z^*))
 \]
 takes the form
  \[
 \lambda^{\K}+\lambda^{\Omega}=0, \ \ (\lambda^{\K}, \lambda^{\Omega})\in
 -\K\times N_{\Omega}(0,0).
 \]
 Take $\lambda^{\K}=(-1,0)$ and $\lambda^G=(1,0)$. Then $\lambda^{\K}\in -\K$ and $\lambda^{\Omega}\in N_{\Omega}(0,0)$. Hence NNAMCQ  does not hold at $z^*$.

 It is interesting to note that each of the set-valued mappings for the two split systems $M_1(z):=(z_1,z_2)-\K$ and $M_2(z):=(z_1,z_2)-\Omega$ are both metrically regular at $(z^*,0)$, but the one for the whole system $M(z)=(z,z)-\K\times \Omega$ is only metrically subregular  (not metrically regular) at $(z^*,0)$.
 \end{example}

In many situations, the constraint system $P(z)\in D$ can be split into subsystems $P_1(z)\in D_1, P_2(z)\in D_2$  such that one subsystem can be checked to have error bound property easily. In Klatte and Kummer \cite[Theorem 2.5]{KK02}, it is shown that if both $M^{-1}_1$ and $M^{-1}_2$ are calm at $(0, z^*)$ and  $M_2$ is pseudo-Lipschitz continuous at $(z^*,0)
$,  then checking the calmness of the intersection $M^{-1}(\alpha, \beta):=M^{-1}_1(\alpha)\cap M^{-1}_2(\beta)$  at $(0, 0,z^*)$ can be reduced to checking the calmness of $H(\beta):=M^{-1}_1(0)\cap M^{-1}_2(\beta)$  at $(0, z^*)$. In Example \ref{ex4.1}, both $M^{-1}_1$ and $M_2^{-1}$ are calm at $(0, z^*)$ and $M_2$ can be checked to be pseudo-Lipschitz continuous at $(z^*,0)$ by using Mordukhovich criterion,
 and $H(\beta)=\{z|z\in {\cal K}, z-\beta\in \Omega\}$ is calm at $(0,z^*)$ as a polyhedral multifunction. This ensures the calmness of $M^{-1}$, or equivalently, the metric subregularity of the whole system $M(z)=P(z)-D$.
 %In addition, the FOSCMS provides the verifiable conditions via the initial problem data to study the metric subregularity for non-polyhedral multifunctions.

 In \cite[Theorem 2]{GfrJYe},  the first order sufficient condition for metric subregularity for a split system with product of two sets is given. When one of the subsystem is known to be metrically subregular, the condition given in \cite[Theorem 2]{GfrJYe} is completely verifiable using the initial data of the problem.
We now extend this result  to the product of finitely many sets.
\begin{thm}\label{PropEquGrSubReg}
  {Let  $P(z^*)\in D$
  %Let $\zb\in\Omega$ be feasible for \eqref{EqGenOptProbl}
  and assume that $P$ and $D$} can be written in the form
  \[P(z)=(P_1(z),P_2(z),\dots, P_I(z)), \quad   D=D_1\times D_2\times \dots \times D_I,\]
where $P_i:\R^n\to\R^{s_i}$ are smooth and $D_i\subseteq \R^{s_i}$, $i=1,2,\dots, I$, are closed  such that the set-valued map $M_1(z):=P_1(z)-D_1$ is metrically subregular at $(z^*,0)$. Further assume that for every $0\not=w$ such that $\nabla P_i(z^*)w\in T_{D_i}(P_i(z^*))$, $i=1,2,\dots, I$, one has
     \begin{eqnarray}
  &&   \left. \begin{array}{l}  \nabla P_1(z^*)^T\lambda^1+\displaystyle \sum_{i=2}^I \nabla P_i(z^*)^T\lambda^i=0,\;\\
  \lambda^i\in N_{D_i}(P_i(z^*);\nabla P_i(z^*)w)\, \ \forall i=1,2,\dots,  I \end{array} \right\}  \Longrightarrow\;\lambda^i=0 \ \ \forall  i=2,\dots, I. \label{eqn-add-1}
  \end{eqnarray}
      Then the set-valued mapping $M(z):=P(z)-D$ is metrically subregular at $(z^*,0)$. Moreover, if all $D_i$ except at most  $D_1$ are directionally regular at $P_i(z^*)$, then  (\ref{eqn-add-1}) is equivalent to that
       for every $w\neq 0$ such that $\nabla P(z^*)w\in T_{D}(P(z^*))$, one has
     \begin{eqnarray}
&& \left \{ \begin{array}{l}
     \nabla P_1(z^*)^T\lambda^1+\displaystyle \sum_{i=2}^I \nabla P_i(z^*)^T\lambda^i=0,\;\\
   \lambda^1\in N_{D_1}(P_1(z^*);\nabla P_1(z^*)w), \\
 ( \lambda^2,\dots, \lambda^I)\in N_{D_2\times \dots \times {D_I}}\big((P_2(z^*), \dots, P_I(z^*));\nabla P_2(z^*)w,\dots, \nabla P_I(z^*)w\big)
 \end{array} \right. \label{Thm4.4}
 \\
  && \Longrightarrow\;\lambda^i=0 \ \forall  i=2,\dots, I. \nonumber\end{eqnarray}
\end{thm}
 \beginproof
 Let $w\not =0$ satisfying $\nabla P(z^*)w\in T_{D}(P(z^*))$
such that  (\ref{Thm4.4}) holds. Then by Proposition \ref{TanNor}, we have $\nabla P_i(z^*)w\in T_{D_i}(P_i(z^*)), $ and $ \lambda^i\in N_{D_i}(P_i(z^*);\nabla P_i(z^*)w)$ for $i=1,2,\dots, I$. Since (\ref{eqn-add-1}) holds at $z^*$, it follows that $\lambda^i=0$ for   $i=2,\dots, I.$
Applying
 \cite[Theorem 2]{GfrJYe}, we have that the set-valued mapping $M(z):=P(z)-D$ is metrically subregular at $(z^*,0)$.  Moreover suppose that all $D_i$ except at most  $D_1$ are directionally regular at $P_i(z^*)$. Then by Proposition \ref{TanNor},
 \begin{eqnarray*}
&& T_D(P(z^*))= T_{D_1}(P_1(z^*))\times \dots \times T_{D_I}(P_I(z^*)),\\
&& { N_{D_2\times \dots \times {D_I}}\big((P_2(z^*), \dots, P_I(z^*));\nabla P_2(z^*)w,\dots, \nabla P_I(z^*)w\big)}\\
&&=N_{D_2}(P_2(z^*);\nabla P_2(z^*)w)\times \dots \times N_{D_I}(P_I(z^*);\nabla P_I(z^*)w),
\end{eqnarray*}
and hence the two conditions are equivalent.
 \endproof

%is \cite[Lemma 1]{GfrKl16}.

\section{Expressions for tangent cones}
In order to use the sufficient conditions for metric subregularity in terms of directional normal cones, one needs to derive the formula for the tangent cone involved. In this section we derive the exact expressions for the tangent cone of the second-order cone complementarity set. Moreover we show that it is geometrically derivable.

The following formula for the tangent cone of the second-order cone is well-known.
\begin{prop}\cite[Lemma 25]{BR05}\label{tcone} Let ${\cal K}$ be the $m$-dimensional second-order cone.
\begin{eqnarray*}
 T_{\cal K}(x)
 &=& \left\{
 \begin{array}{lll}
 & \mathbb{R}^m \ &{\rm if } \  x\in {\rm int}{\cal K}  \\
 & {\cal K} \ &{\rm if}\ x=0 \\
 &  d\in \mathbb{R}^m: -d_1+\bar{x}^T_2 d_2\leq 0 \ \ & {\rm if}\ x\in {\rm bd}{\cal K}\setminus\{0\}
 \end{array}
\right\}.
 \end{eqnarray*}
\end{prop}

%According to
% Theorem  \ref{PropEquGrSubReg}, we need to calculate the tangent cone to the set $\Omega$, defined as
Let ${\cal K}$ be the $m$-dimensional second-order cone and   \begin{eqnarray} \label{Omega}
 \Omega = \{(x,y) \, | \, \K \ni x \perp y \in \K\}
 \end{eqnarray}
 be the corresponding second-order cone complementarity set.
In what follows we show that the set $\Omega$ is geometrically derivable and give a characterization in terms of the metric projection operator. The characterization of the tangent cone was also given in \cite[Proposition 3.1]{JLZ15}.
 \begin{prop}\label{prop-tangent} The set $\Omega$ is geometrically derivable, and for any $(x,y)\in \Omega$,
 \[
 T_{\Omega}(x,y)=T^i_\Omega(x,y)= \left\{ (d,w) \, | \, \Pi'_{\cal K}(x-y;d-w)=d \right\}.
 \]
 \end{prop}

 \beginproof Since $T^i_{\Omega}(x,y)\subseteq T_\Omega(x,y)$, it suffices to prove $$T_\Omega(x,y)\subseteq\Upsilon(x,y) \subseteq T^i_{\Omega}(x,y),$$
 where $\Upsilon(x,y):=\left \{ (d,w) \, | \, \Pi'_{\cal K}(x-y;d-w)=d \right \}$.
 Take $(d,w)\in T_\Omega(x,y)$. Then by definition, there exist  $t_n\downarrow 0$ and
 $(d_n,w_n)\to (d,w)$ such that
 $(x,y)+t_n (d_n,w_n)\in \Omega$. By (\ref{relation-1}), we have
 \[
 \Pi_{\cal K}(x+t_n d_n - y-t_n w_n)=x+t_n d_n=\Pi_{\cal K}(x-y)+t_n d_n.
 \]
 Hence
 \[
 \Pi_{\cal K}'(x-y;d-w)=\lim\limits_{{t_n}\downarrow 0}
 \frac{\Pi_{\cal K}(x+t_n d_n - y-t_n w_n)-\Pi_{\cal K}(x-y)}{t_n}=d.
 \]
Therefore
 \(
 T_\Omega(x,y) \subseteq\Upsilon(x,y).
 \)
 Now take $(d,w)\in \Upsilon(x,y)$.
 Then for any given $t_n\downarrow 0$
 \[
 \Pi_{\cal K}(x-y+t_n(d-w))-\Pi_{\cal K}(x-y)
 = t_n \Pi_\K'(x-y;d-w)+o(t_n)
 =t_n d+o(t_n),
 \]
 i.e.,
 \[
 \Pi_{\cal K}(x-y+t_n(d-w))=\Pi_{\cal K}(x-y)+t_n d+o(t_n).
 \]
 Hence
 \begin{eqnarray*}
 & & \bigg( \Pi_{\cal K}(x-y)+t_n d+o(t_n), \Pi_{\cal K}(x-y)+t_n d + o(t_n)
     -(x-y+t_n(d-w)) \bigg)\in \Omega \\
 & \Longleftrightarrow &
 \bigg( x+t_n d+o(t_n), y +t_n w+o(t_n) \bigg)\in \Omega.
 \end{eqnarray*}
 This means $(d,w)\in \liminf\limits_{t_n\downarrow 0} \ds \frac{\Omega-(x,y)}{t_n}$. Hence
 $\Upsilon(x,y)\subseteq T^i_\Omega(x,y)$.
 \endproof

  With the above result, an explicit expression of tangent cone $T_\Omega$
 is also given in \cite[Theorem 3.1]{JLZ15}. However, in that explicit
 formula, for the case where $x-y\notin \K\cup \K^\circ$, the directional derivative of
 the projection operator is involved. In the hope that only the initial
 data on $x,y$ is used, we next provide another explicit expression for
 $T_\Omega(x,y)$ without involving $\Pi_{\K}'(x-y)$.
 \medskip

%---------------------------------------------------------------------------------- Theorem 2.1
 \begin{thm}\label{formula-regular normal cone}
 Let $\Omega$ be defined as in (\ref{Omega}). Then for any $(x,y)\in \Omega$,
 %the structure of the tangent cone of $\Omega$ at $(x,y)$ is described as
 \begin{eqnarray*}
 \lefteqn{T^i_\Omega(x,y)=T_\Omega(x,y)}\\
 &=& \left\{ (d,w) \left|
 \begin{array}{lll}
 & d=0,\ w\in \mathbb{R}^m, \ &{\rm if} \ x=0, \ y\in {\rm int}{\cal K};  \\
 & d\in \mathbb{R}^m,\ w=0, \ &{\rm if}\ x\in {\rm int}{\cal K}, \ y=0; \\
 & x_1\hat{w}-y_1d\in \mathbb{R} x, \ \ d\perp y,\ w\perp x, \ \ & {\rm if}
   \ x,y \in {\rm bd}{\cal K}\backslash\{0\}; \\
 & d=0,\ w\in T_{\cal K}(y) \ {\rm or}\ d\in \mathbb{R}_{+}\hat{y}, \ w\perp \hat{y}, \ \
 & {\rm if}\ x=0,\ y\in {\rm bd}{\cal K}\backslash \{0\}; \\
 & d\in T_{\cal K}(x),\ w=0 \ {\rm or}\ d\perp \hat{x}\ , w\in \mathbb{R}_{+}\hat{x}, \ \
 & {\rm if}\ x\in {\rm bd}{\cal K}\backslash \{0\},\ y=0; \\
 & d\in {\cal K},\ w\in {\cal K}, \ d\perp w, \ \ & {\rm if}\ x=0, \ y=0.
 \end{array}
 \right.\right\}.
 \end{eqnarray*}
 \end{thm}

 \beginproof
 Note that $\Pi_{\cal K}(x-y)$ is continuously differentiable at $x-y$, {provided that}
 $x=0$ and $y\in {\rm int}{\cal K}$, or $x\in {\rm int}{\cal K}$ and
 $y=0$, or $x,y\in {\rm bd}{\cal K}\backslash \{0\}$ with $x^T y=0$.
 Hence $D^*\Pi_{\cal K}(x-y)=\nabla \Pi_{\cal K}(x-y)$ in the
 above cases, which in turn implies
 \[
 \Pi'_{\cal K}(x-y;d-w)=d
 \ \Longleftrightarrow \
 D^*\Pi_{\cal K}(x-y)(d-w)=d
 \ \Longleftrightarrow \
 (w,-d)\in N_\Omega (x,y),
 \]
 where the second equivalence is due to Proposition 2.4. By the
 expression of the limiting normal cone in Proposition \ref{formula-normal cone}, we have the following conclusions.

 \medskip
 \noindent
 \underline{Case 1}: If $x=0$ and $y\in {\rm int}{\cal K}$, then
 $w\in \mathbb{R}^m$ and $d=0$.

 \medskip
 \noindent
 \underline{Case 2}: If $x\in {\rm int}{\cal K}$ and $y=0$, then
 $d\in \mathbb{R}^m$ and $w=0$.

 \medskip
 \noindent
 \underline{Case 3}: If $x,y\in {\rm bd}{\cal K}\backslash \{0\}$ with
 $x^T y=0$, then
 \[
 x_1\hat{w}-y_1d\in \mathbb{R} x, \ \ d\perp y,\ w\perp x.
 \]

 \noindent
 \underline{Case 4}: If $x=0$ and $y\in {\rm bd}{\cal K}\backslash \{0\}$,
 since by Proposition \ref{prop-tangent}, $(d,w)\in T_\Omega(x,y)$ is equivalent to saying that
 $\Pi'_{\cal K}(x-y;d-w)=d$. According to the  formula of
 directional derivative for $\Pi_{\cal K}$ in Proposition \ref{directd}(iii), we have $(d,w)\in T_\Omega(x,y)$ if and only if
 \begin{equation}\label{tangent-case-4}
 \frac{1}{2}\bigg( d_1-w_1-\bar{y}_2^Td_2+\bar{y}_2^Tw_2  \bigg)_+
 \left[
 \begin{array}{c}
 1 \\
 -\bar{y}_2
 \end{array}
 \right]
 =\left[
 \begin{array}{c}
 d_1 \\
 d_2
 \end{array}
 \right].
 \end{equation}
We now claim that the set of solutions to equation (\ref{tangent-case-4}) is
 \begin{equation} \label{equivalent-solution-tangent}
 \Big \{ (d,w) \, | \, d=0,\ w\in T_{\cal K}(y) \ {\rm or}\
 d \in \mathbb{R}_{+}\hat{y}, \ w \perp \hat{y} \Big \}.
 \end{equation}
 \medskip
 \noindent
 By definition of the tangent cone, $(d,w)\in T_\Omega (x,y)$ if and only if there exists $(d_n,w_n)\to (d,w)$
 and $t_n\to 0$ with $t_n\geq 0$ satisfying
 $(t_n d_n, y+t_n w_n)=(x,y)+t_n(d_n,w_n)\in \Omega$, i.e.,
 \[
 t_n d_n \in {\cal K}, \quad y+t_n w_n\in {\cal K}, \quad  t_n d_n \perp y+t_n w_n,
 \]
 which implies
 \[
 d\in {\cal K}, \quad w\in T_{\cal K}(y), \quad d\perp y.
 \]
% If $(d,w)\in T_\Omega (x,y)$, from the above arguments, then $(d,w)$
% satisfies (\ref{def-tangent-cone}).
Note that  $d\in {\cal K}$ and
 $\langle y,d \rangle=0$ implies that either  $d=0$ or $d\in \mathbb{R}_{++}\hat{y}$.
 If $d=0$, then (\ref{tangent-case-4}) takes the form as
 \[
 \frac{1}{2} \bigg( -w_1+\bar{y}_2^Tw_2  \bigg)_+
 \left[
 \begin{array}{c}
 1 \\
 -\bar{y}_2
 \end{array}
 \right]
 =\left[
 \begin{array}{c}
 0 \\
 0
 \end{array}
 \right],
 \]
 which implies $-w_1+\bar{y}_2^Tw_2\leq 0$, i.e., $w\in T_{\cal K}(y)$
 by the formula of tangent cone of ${\cal K}$ in Proposition \ref{tcone}.
 If $d\in \mathbb{R}_{++}\hat{y}$, then $d=\tau \hat{y}=\tau (y_1,-y_2)$ for some $\tau>0$. Hence
 $$\bar{y}_2^Td_2=\bar{y}_2^T(-\tau y_2)=-\tau \|y_2\|=-\tau y_1=-d_1.$$
 It then follows
 from (\ref{tangent-case-4}) that
 \[
 d_1
 = \frac{1}{2} \bigg( d_1-w_1-\bar{y}_2^Td_2+\bar{y}_2^Tw_2  \bigg)_+
 = \frac{1}{2} \bigg( 2d_1-w_1+\bar{y}_2^Tw_2  \bigg)_+ , \]
 which implies $-w_1+\bar{y}_2^Tw_2=0$, i.e., $w\perp \hat{y}$.
 Thus, $(d,w)$ satisfies (\ref{equivalent-solution-tangent}).

 Conversely,
 if $d=0$ and $w\in T_{\cal K}(y)$, then $-w_1+\bar{y}_2^Tw_2\leq 0$ by Proposition \ref{tcone}, and hence
 \[(d_1-w_1-\bar{y}_2^Td_2+\bar{y}_2^Tw_2)_+=(-w_1+\bar{y}_2^Tw_2)_+= 0,\]
 which implies that
 (\ref{tangent-case-4}) holds, i.e., $(d,w)\in T_\Omega(x,y)$. For the
 other case, if $d\in \mathbb{R}_{++}\hat{y}$ and $w\perp \hat{y}$, then
 $d=\tau \hat{y}$ for some $\tau >0$ and $w_1=\bar{y}_2^T w_2$. Therefore
 \begin{eqnarray*}
 \frac{1}{2} \bigg( d_1-w_1-\bar{y}_2^Td_2+\bar{y}_2^Tw_2  \bigg)_+
 \left[
 \begin{array}{c}
 1 \\
 -\bar{y}_2
 \end{array}
 \right]
 &=& \frac{1}{2} \bigg( d_1-\bar{y}_2^Td_2  \bigg)_+
 \left[
 \begin{array}{c}
 1 \\
 -\bar{y}_2
 \end{array}
 \right] \\
 &=& d_1 \left[
 \begin{array}{c}
 1 \\
 -\bar{y}_2
 \end{array}
 \right] \\
 &=& \left[
 \begin{array}{c}
 d_1 \\ d_2
 \end{array}
 \right],
 \end{eqnarray*}
 where the second equation is due to
 $d_1-\bar{y}_2^Td_2=d_1+\tau\bar{y}_2^Ty_2=d_1+\tau\|y_2\|=d_1+\tau y_1=2d_1$
 and the third equation comes from
 $-d_1 \bar{y}_2=-\tau y_1 \bar{y}_2=-\tau y_2=d_2$. This means that (\ref{tangent-case-4})
 holds, i.e., $(d,w)\in T_\Omega(x,y)$. In summary, we have shown that
 \[
 T_\Omega(x,y)
 = \big\{ (d,w) \, \big| \, d=0,\ w\in T_{\cal K}(y) \ {\rm or} \
 d \in \mathbb{R}_{+}\hat{y}, \ w\perp \hat{y} \big\}.
 \]

 \medskip
 \noindent
 \underline{Case 5}: If $x\in {\rm bd}{\cal K}\backslash \{0\}$ and $y=0$,
 by symmetry to Case 4, we have
 \[
 T_\Omega(x,y)=\big \{ (d,w) \, \big | \, d \in T_{\cal K}(x),\ w=0 \ {\rm or} \
 d \perp \hat{x},\ w\in \mathbb{R}_{+}\hat{x} \big \}.
 \]

 \medskip
 \noindent
 \underline{Case 6}: If $x=0$ and $y=0$, then
 $\Pi_{\cal K}'(0;h)=\Pi_{\cal K}(h)$ by Proposition \ref{directd}(iv). It follows from Proposition \ref{prop-tangent} that
 \[
 (d,w)\in T_\Omega(x,y) \ \ \Longleftrightarrow \ \ \Pi_{\cal K}(d-w)=d
 \ \ \Longleftrightarrow \ \ -w\in N_{\cal K}(d)
 \ \ \Longleftrightarrow \ \ (d,w)\in \Omega,
 \]
 i.e., $T_\Omega(x,y)=\Omega$. In fact, this case can also be obtained by noting
 that $\Omega$ is a cone.
 \endproof

 \medskip

When $m=1,2$, the tangent cone $T_\Omega$ have simpler expression given below.
{For example, when $m=1$, the second-order cone complementarity set $\Omega$ is reduced to the vector complementarity set
$\{(a,b)\in \mathbb{R}^2| a\geq 0,b\geq 0, ab=0\}$, and hence ${\rm bd}{\K}\backslash \{0\}$ is empty; when $m=2$, the condition
$x_1\hat{w}-y_1d\in \mathbb{R} x$ can be dropped.}

 \medskip

%-------------------------------------------------------------------------------- Corollary 2.1
 \begin{cor}\label{co-tangent-cone}
 Let $\Omega$ be defined as in (\ref{Omega}). If $m=1$, then
 \[
 T_\Omega(x,y)=\left\{(d,w)\left| \begin{array}{lll}
  &d=0,\  \ &{\rm if} \ x=0, \  y>0 \\
 & w=0, \ & {\rm if}\ x>0, \ y=0 \\
 & d\geq 0,\ w\geq 0, \ d\perp w, \ \ & {\rm if}\ x=0,  \ y=0
 \end{array}\right.\right\}.
 \]
 If $m=2$, then
 \[
 T_\Omega(x,y)=\left\{(d,w)\left|\begin{array}{lll}
 & d=0,\ w\in \mathbb{R}^2, \ &{\rm if} \ x=0,\ y\in {\rm int}{\cal K} \\
 & d\in \mathbb{R}^2,\ w=0, \ &{\rm if}\ x\in {\rm int}{\cal K}, \ y=0 \\
 & d\perp y,\ w\perp x, \ \ & {\rm if} \ x, y\in {\rm bd}{\cal K}\backslash\{0\}
  \\
 & d=0,\ w\in T_{\cal K}(y) \ {\rm or}\ d\in \mathbb{R}_{+}\hat{y}, \ w\perp \hat{y},
   \ \ &{\rm if}\ x=0,\  y\in {\rm bd}{\cal K}\backslash \{0\} \\
 & d\in T_{\cal K}(x),\ w=0 \ {\rm or}\ d\perp \hat{x}\ , w\in \mathbb{R}_{+}\hat{x},
   \ \ &{\rm if}\ x\in {\rm bd}{\cal K}\backslash \{0\},\  y=0 \\
 & d\in {\cal K},\ w\in {\cal K}, \ d\perp w, \ \ &{\rm if}\ x=0, \  y=0
 \end{array}
 \right.\right\}.
 \]
 \end{cor}

 \beginproof
 If $m=1$, then $\K=\mathbb{R}_+$, and hence ${\rm bd}\K\backslash \{0\}=\emptyset$.
 Thus the desired result follows.

 \medskip
 \noindent
 If $m=2$, we show that the condition $x_1\hat{w}-y_1d\in \mathbb{R} x$
 is implied by $d\perp y$ and $w\perp x$ in the case of
 $x,y\in {\rm bd}\K\backslash \{0\}$ and $x^Ty=0$. In fact, if $x_1=x_2>0$, then $y_1=-y_2$,
 and hence $w_1+w_2=0$ and $d_1-d_2=0$. Thus
 \[
 x_1 \left[
 \begin{matrix}
 w_1 \\
 -w_2
 \end{matrix}
 \right]
 -y_1 \left[
 \begin{matrix}
 d_1 \\
 d_2
 \end{matrix}
 \right]
 = x_1 \left[
 \begin{matrix}
 w_1 \\
 w_1
 \end{matrix}
 \right]
 -y_1 \left[
 \begin{matrix}
 d_1 \\
 d_1
 \end{matrix}
 \right]
 = \frac{x_1w_1-y_1d_1}{x_1}
 \left[
 \begin{matrix}
 x_1 \\
 x_2
 \end{matrix}
 \right].
 \]
 Similarly, if $x_1=-x_2$, then $y_1=y_2$, and hence $d_1+d_2=0$ and $w_1-w_2=0$.
 Thus
 \[
 x_1 \left[
 \begin{matrix}
 w_1 \\
 -w_2
 \end{matrix}
 \right]
 -y_1 \left[
 \begin{matrix}
 d_1 \\
 d_2
 \end{matrix}
 \right]
 =x_1 \left[
 \begin{matrix}
 w_1 \\
 -w_1
 \end{matrix}
 \right]
 -y_1 \left[
 \begin{matrix}
 d_1 \\
 -d_1
 \end{matrix}
 \right]
 = \frac{x_1w_1-y_1d_1}{x_1}
 \left[
 \begin{matrix}
 x_1 \\
 -x_1
 \end{matrix}
 \right]
 = \frac{x_1w_1-y_1d_1}{x_1}
 \left[
 \begin{matrix}
 x_1 \\
 x_2
 \end{matrix}
 \right].
 \]
 \endproof

 \section{Expressions for directional   normal cones}
 In order to use  FOSCMS for the second-order cone complementarity system, one needs to derive the {exact formula} for the directional normal cone of the second-order cone complementarity set. Moreover these results are of their own interest.

By formulating the vector complementarity set as the union of finitely many polyhedral convex sets, the formula of the directional normal cone of the vector complementarity set is given in \cite[Lemma 4.1]{Gfr14}. In contrast to the vector complementarity set, the second-order cone complementarity set  cannot be represented as the  union of finitely many polyhedral convex sets.
 In the following theorem we derive an {explicit} expression for the directional normal cone for the $m$-dimensional second-order cone complementarity set $\Omega$ defined as in (\ref{Omega}). Note that in the case where $m=1$,  ${\rm bd}\K\backslash \{0\}=\emptyset$ and hence the formula we derived reduced to the one given in \cite[Lemma 4.1]{Gfr14} for this case.
 %With the non-polyhedral structure of $\K$, the analysis .
 %Recall that the formula of the tangent cone and the normal cone of $\Omega$ has been established in Theorem \ref{formula-regular normal cone} and Proposition \ref{formula-normal cone} respectively.

 \begin{thm}\label{formula-directional normal} The second-order cone complementarity set is directionally regular. For any $(x,y)\in \Omega$ and $(d,w)\in T_{\Omega}(x,y)=T_{\Omega}^i(x,y)$, the directional  normal cone can be calculated as follows.

\noindent Case 1: $x=0,y\in {\rm int}\K$,
\[
N_{\Omega}\big((x,y);(d,w)\big)=N_{\Omega}(x,y)= \mathbb{R}^m\times \{0\}.
\]
\noindent Case 2: $x\in {\rm int}\K,y=0$,
\[
N_{\Omega}\big((x,y);(d,w)\big)=N_{\Omega}(x,y)=\{0\} \times \mathbb{R}^m.
\]
\noindent Case 3: $x,y\in {\rm bd}\K\backslash\{0\},$
\[
N_{\Omega}\big((x,y);(d,w)\big)=N_{\Omega}(x,y).
\]
\noindent Case 4:  $x\in {\rm bd}\K\backslash \{0\},y=0$,
\begin{eqnarray*}
 N_{\Omega}\big((x,y);(d,w)\big)=\left\{(u,v)\left|
 \begin{array}{lll}
  u=0, v\in \mathbb{R}^m \ &{\rm if} \ d\in {\rm int}T_{\K}(x),w=0 \\
  N_{\Omega}(x,y) \ &{\rm if}\ d\in {\rm bd}T_{\K}(x), w=0 \\
   u\in \mathbb{R} \hat{x},\
 v\perp \hat{x}\ & {\rm if}\ d\perp \hat{x}, w\in \mathbb{R}_+\hat{x}\backslash \{0\}
 \end{array}
 \right.\right\}.
 \end{eqnarray*}
\noindent Case 5: $x=0,y\in {\rm bd}\K\backslash \{0\},$
\begin{eqnarray*}
 N_{\Omega}\big((x,y);(d,w)\big)=\left\{(u,v)\left|
 \begin{array}{lll}
  u\in \mathbb{R}^m, v=0 \ &{\rm if} \ d=0, w\in {\rm int}T_{\K}(x) \\
  N_{\Omega}(x,y) \ &{\rm if}\ d=0, w\in {\rm bd}T_{\K}(x) \\
   v\in \mathbb{R} \hat{y},\
 u\perp \hat{y}\ & {\rm if}\ d\in \mathbb{R}_+\hat{y}\backslash \{0\}, w\perp \hat{y}
 \end{array}
 \right.\right\}.
 \end{eqnarray*}
\noindent Case 6:
 $x=0,y=0$,
 \[
 N_{\Omega}\big((x,y);(d,w)\big)=N_{\Omega}(d,w).
 \]
 Here the formula of the tangent cone and the normal cone of $\Omega$  are given as in Theorem \ref{formula-regular normal cone} and Proposition \ref{formula-normal cone} respectively.
\end{thm}

\beginproof
In Cases 1-3, since  it always has\[
N_{\Omega}^i\big((x,y);(d,w)\big) \subseteq N_{\Omega}\big((x,y);(d,w)\big)\subseteq N_{\Omega}(x,y),
\]
it suffices to show that
\begin{equation}
 N_{\Omega}(x,y) \subseteq N_{\Omega}^i\big((x,y);(d,w)\big).\label{inclusion}
\end{equation} For any $(u,v)\in N_{\Omega}(x,y)$, in order to show that $(u,v) \in
N^i_{\Omega}\big((x,y);(d,w)\big)$,  for  any sequences $t_n\downarrow 0$, we need to find   $(d^n,w^n)\to (d,w)$ and $(u^n,v^n)\to (u,v)$ satisfying $(u^n,v^n)\in \widehat{N}_{\Omega}(x+t_nd^n,y+t_nw^n)$.

\noindent {\bf Case 1.} $x=0$, $y\in {\rm int}\K$.
%Note that \[N_{\Omega}((x,y);(d,w))\subset N_{\Omega}(x,y)=\{(u,v)| v=0\}.\]
Since $(u,v) \in  N_{\Omega}(x,y), (d,w)\in T_{\Omega}(x,y)$, then by Theorem \ref{formula-regular normal cone}, $d=0, w\in \mathbb{R}^m$ and by Proposition \ref{formula-normal cone}, $u\in \mathbb{R}^m, v=0$.  By letting
\[(u^n,v^n):=(u,v)=(u,0) \ \ {\rm and} \ \ (d^n,w^n):=(d,w)=(0,w),\]
 we have $y+t_nw^n\in {\rm int}\K$ for $n$ sufficiently large, and hence \[(u^n,v^n)=(u,0)\in \widehat{N}_{\Omega}(0,y+t_nw)= \widehat{N}_{\Omega}(x+t_nd^n,y+t_nw^n).\] Hence (\ref{inclusion}) holds.
 %So $(u,v)\in N_{\Omega}((x,y);(d,w))$.

\noindent {\bf Case 2.} $x\in {\rm int}\K$, $y=0$. This case is symmetric to Case 1 and we omit the proof.
%\[N_{\Omega}((x,y);(d,w))= \{(u,v)| u=0\}=N_{\Omega}(x,y).\]

\noindent {\bf Case 3.} $x,y\in {\rm bd}\K\backslash\{0\}$. Then $x-y\in (-\K\cup \K)^c$. Since $(d,w)\in T_\Omega(x,y)=T_\Omega^i(x,y)$, by definition of the inner tangent cone,
for any $t_n\downarrow 0$, there exists $(d^n,w^n)\to (d,w)$ such that $(x,y)+t_n(d^n,w^n)\in \Omega$. We now construct a sequence $(u^n,v^n)$ such that
$(u^n,v^n)\to (u,v)$ and $(u^n,v^n)\in \widehat{N}_{\Omega}(x+t_nd^n,y+t_nw^n)$.  By Proposition \ref{Lem3.1},
$y=k\hat{x}$ with $k=y_1/x_1$. Hence $x_1-y_1=(1-k)x_1, x_2-y_2=(1+k)x_2$. By Proposition \ref{projection}(iii), the metric projection $\Pi_{\K}$ is differentiable at $x-y$ and
\[\nabla \Pi_{\K}(x-y)=\frac{1}{2}\begin{bmatrix}
1 & \bar{x}_2^T \\
\bar{x}_2 & I+\frac{1-k}{1+k}(I-\bar{x}_2\bar{x}_2^T)
\end{bmatrix}=\frac{1}{2}\left\{\begin{bmatrix}
1 & \bar{x}_2^T \\
\bar{x}_2^T & I
\end{bmatrix}+\frac{1-k}{1+k}\begin{bmatrix}
0 & 0 \\
0 & I-\bar{x}_2\bar{x}_2^T
\end{bmatrix}\right\}.\]
By
\cite[Lemma 1]{CCT04}, the eigenvalue values of the matrix  are $0,1$ and $1/1+k$ with multiplicity $n-2$ .

\underline{Case 3(i)}: If $k\neq 1$, then the eigenvalue of the matrix $I-2\nabla \Pi_{\K}(x-y)$ is $1, -1, \frac{k-1}{k+1}$. So $I-2\nabla \Pi_{\K}(x-y)$ is invertible. Since $\Pi_{\K}$ is continuously differentiable at $x-y$,  $I-2\nabla \Pi_{\K}(x+t_nd^n-y-t_nw^n)$  is  also invertible for sufficiently large $n$.
Let
\[\alpha(n):=\bigg(I-2\nabla \Pi_{\K}(x+t_nd^n-y-t_nw^n)\bigg)^{-1}\bigg(\nabla \Pi_{\K}(x+t_nd^n-y-t_nw^n)-\nabla \Pi_{\K}(x-y)\bigg)(-u-v).\]
Then
$\alpha(n)\to 0$ as $n\to \infty$. Note that
\begin{eqnarray*}
\bigg(I-2\nabla \Pi_{\K}(x+t_nd^n-y-t_nw^n)\bigg)\alpha(n)=\bigg(\nabla \Pi_{\K}(x+t_nd^n-y-t_nw^n)-\nabla \Pi_{\K}(x-y)\bigg)(-u-v)
\end{eqnarray*}
i.e.,
\begin{eqnarray*}
&&\alpha(n)\\
&=&\bigg(\nabla \Pi_{\K}(x+t_nd^n-y-t_nw^n)-\nabla \Pi_{\K}(x-y)\bigg)(-u-v)+2\nabla \Pi_{\K}(x+t_nd^n-y-t_nw^n)\alpha(n)\\
&=& \nabla \Pi_{\K}(x+t_nd^n-y-t_nw^n)(-u-v+2\alpha(n))-\nabla \Pi_{\K}(x-y)(-u-v)\\
&=& \nabla \Pi_{\K}(x+t_nd^n-y-t_nw^n)(-u-v+2\alpha(n))-(-v),
\end{eqnarray*}
where the last step is due to $-v=\nabla \Pi_{\K}(x-y)(-u-v)$ since $(u,v)\in N_{\Omega}(x,y)$ by Proposition \ref{Prop2.2}. Hence
\begin{eqnarray*}
-v+\alpha(n)
= \nabla \Pi_{\K}(x+t_nd^n-y-t_nw^n)(-u+\alpha(n)-v+\alpha(n)).
\end{eqnarray*}
Let $(u^n,v^n):=(u-\alpha(n),v-\alpha(n))$. Then
\[(u^n,v^n)\to (u,v) \ \ {\rm and} \ \  -v^n\in \widehat{D}^*\Pi_{\K}(x+t_nd^n-y-t_nw^n)(-u^n-v^n).\] 
 By  Proposition \ref{Prop2.2} $(u^n,v^n)\in \widehat{N}_{\Omega}(x+t_nd^n,y+t_nw^n)$. Hence (\ref{inclusion}) holds.

\underline{Case 3(ii)}: If $k=1$, then the eigenvalue of $I-2\nabla \Pi_{\K}(x-y)$ is $1, -1, \frac{k-1}{k+1}=0$ and hence the matrix $I-2\nabla \Pi_{\K}(x-y)$ is not invertible and the construction of $(u^n,v^n)$ in case 3(i) fails. Note that in this case the eigenvalue of the matrix $I-3\nabla \Pi_{\K}(x-y)$ is $1, -2, -\frac{1}{2}$. So $I-3\nabla \Pi_{\K}(x-y)$ has inverse. We then construct the sequence by taking
$(u^n,v^n):=(u-2\alpha(n),v-\alpha(n))$ with
\[\alpha(n):=\bigg(I-3\nabla \Pi_{\K}(x+t_nd^n-y-t_nw^n)\bigg)^{-1}\bigg(\nabla \Pi_{\K}(x+t_nd^n-y-t_nw^n)-\nabla \Pi_{\K}(x-y)\bigg)(-u-v),\]
and the desired result follows similarly.

\noindent {\bf Case 4.} $x\in {\rm bd}\K\backslash\{0\}$ and $y=0$. Since $(d,w)\in T_{\Omega}(x,y)$, by Theorem \ref{formula-regular normal cone}, there are three possible cases: $w=0, d\in {\rm int}T_{\K}(x)$,  $w=0, d\in {\rm bd}T_{\K}(x)$,  or $d\perp \hat{x}, w\in \mathbb{R}_+\hat{x}\setminus\{0\}$.
%For any sequences $t_n\downarrow 0$, $(d^n,w^n)\to (d,w)$, we have
%\begin{eqnarray*}
%\Limsup\limits_{n\to \infty}\widehat{N}_{\Omega}\big((x+t_nd^n,t_nw^n)\big)
%&=&\Limsup\limits_{{t_nw^n=0}\atop {x+t_nd^n\in {\rm int}\K}}\widehat{N}_{\Omega}\big((x+t_nd^n,t_nw^n)\big)\\
%&&\bigcup
%\Limsup\limits_{{t_nw^n=0}\atop {x+t_nd^n\in {\rm bd}\K\setminus\{0\}}}\widehat{N}_{\Omega}\big((x+t_nd^n,t_nw^n)\big)\\
%&&\bigcup
%\Limsup\limits_{{t_nw^n\in {\rm bd}\K\setminus\{0\}}\atop {x+t_nd^n\in {\rm bd}\K\setminus\{0\}}}\widehat{N}_{\Omega}\big((x+t_nd^n,t_nw^n)\big).
%\end{eqnarray*}

\underline{Subcase 4.1.} $w=0$ and $d\in {\rm int}T_{\K}(x)$.
 Since $(d,w)\in T_{\Omega}(x,y)=T^i_{\Omega}(x,y)$, then for any $t_n\downarrow 0$, there exists $(d^n,w^n)\to (d,w)$ such that $(x,y)+t_n(d^n,w^n)\in \Omega$.
 In this case  $\bar{x}_2^Td_2-d_1<0$. Hence
$\|x_2+t_nd^n_2\|=\|x_2\|+t_n\bar{x}_2^Td^n_2+o(t_n)<x_1+t_nd^n_1$ for  sufficiently large $n$. So $x+t_nd^n\in {\rm int}\K$. It follows that
%Since $(d,w)\in T_\Omega(x,y)$, there exists $t_n\downarrow 0$ and $(d^n,w^n)\to (d,w)$ such that $(x,y)+t_n(d^n,w^n)\in \Omega$. Hence for sufficiently large $n$,
\begin{eqnarray*}
\widehat{N}_{\Omega}\big(x+t_nd^n,t_nw^n\big)
%&=&\Limsup\limits_{{t_nw^n=0}\atop {x+t_nd^n\in {\rm int}\K}}\widehat{N}_{\Omega}\big((x+t_nd^n,t_nw^n)\big)
= \{(u,v)| u=0, v\in \mathbb{R}^m\}.
\end{eqnarray*}
Hence $$N_{\Omega}((x,y);(d,0))=N_{\Omega}^i((x,y);(d,0))= \{(u,v)| u=0, v\in \mathbb{R}^m\}.$$
%Conversely, take $(u,v)$ with $u=0$. Let $(u^n,v^n)=(u,v)$
%and $(d^n,w^n)=(d,w)=(d,0)$. Then $(u^n,v^n)=(u,v)\in \widehat{N}_{\Omega}(x+t_nd,0)=\widehat{N}_{\Omega}(x+t_nd^n,y+t_nw^n)$.
% So $(u,v)\in N_{\Omega}((x,y);(d,w))$. Thus
%\[
%N_{\Omega}\big((x,y);(d,w)\big)=\{(u,v)| u=0, v\in \mathbb{R}^m\}.
%\]

\underline{Subcase 4.2.} $w=0$ and $d\in {\rm bd}T_{\K}(x)$. In this case, it suffices to show
\[
 N_{\Omega}(x,y)\subseteq N_{\Omega}^i\big((x,y);(d,w)\big).
\]
%The $``\subset"$ is clear. Conversely, picking $(u,v)\in N_{\Omega}(x,y)$, we wish to show that $(u,v) \in N_{\Omega}\big((x,y);(d,w)\big)$.
By Proposition \ref{formula-normal cone}, for any $(u,v)\in N_{\Omega}(x,y)$, there are three possible cases: $ u=0, v\in \mathbb{R}^m$ or $ u\in \mathbb{R} \hat{x},\
 v\perp \hat{x}$  or $ u\in \mathbb{R}_-\hat{x},
  \langle v, \hat{x} \rangle \leq 0 $.

\underline{Subcase 4.2(i).} $u=0$ and $v\in \mathbb{R}^m$. Since $d\in {\rm bd}T_{\K}(x)$, we have
$x_1d_1-{x}_2^Td_2=0$ by the formula for the tangent cone. For $t_n\downarrow 0$, let
$\eta(t_n):=\|x_2+t_nd_2\|-\|x_2\|-t_n\bar{x}_2^Td_2$,
$w^n:=0$, {and
\[d^n:=\left\{\begin{array}{ll}
(d_1+2|\eta(t_n)|/t_n,d_2), &\ {\rm if}\ \eta(t_n)\neq 0\\
(d_1+t_n,d_2), &\ {\rm otherwise}.   \end{array}\right.\]
If
$\eta(t_n)\neq 0$, then
\[x_1+t_nd^n_1=
x_1+t_nd_1+2|\eta(t_n)|>\|x_2\|+t_n\bar{x}_2^Td_2+\eta(t_n)=\|x_2+t_nd^n_2\|;\]
otherwise
\[x_1+t_nd^n_1=x_1+t_nd_1+t_n^2=\|x_2\|+t_n\bar{x}_2^Td_2+t_n^2>\|x_2+t_nd_2\|.\]}
Hence $x+t_nd^n\in {\rm int}\K$ and $y+t_nw^n=0$. This ensures $(u,v)\in \widehat{N}_{\Omega}(x+t_nd^n,0) =\widehat{N}_{\Omega}(x+t_nd^n,y+t_nw^n)$. So $(u,v)\in N^i_{\Omega}((x,y);(d,w))$.

\underline{Subcase 4.2(ii).} $u\in \mathbb{R} \hat{x}$ and $v\perp \hat{x}$.  In this case,
$u_2=-u_1\bar{x}_2$ and $v_1-\bar{x}_2^Tv_2=0$. This is equivalent to
\begin{eqnarray*}
\left\{\begin{array}{ll}
u_1=v_1-\bar{x}_2^T(u_2+v_2)\\
u_2=-u_1\bar{x}_2
\end{array} \right. &\Longleftrightarrow & \left\{\begin{array}{ll}
u_1=v_1-\bar{x}_2^T(u_2+v_2)\\
2u_2=-(u_1+v_1)\bar{x}_2+(v_1-u_1)\bar{x}_2
\end{array} \right. \\
&\Longleftrightarrow &  \left\{\begin{array}{ll}
u_1=v_1-\bar{x}_2^T(u_2+v_2)\\
2u_2=-(u_1+v_1)\bar{x}_2+(u_2+v_2)^T\bar{x}_2\bar{x}_2
\end{array} \right. \\
& \Longleftrightarrow &
2u=\begin{bmatrix}
1 & -\bar{x}_2^T\\
-\bar{x}_2 & \bar{x}_2\bar{x}_2^T
\end{bmatrix}(u+v)\\
& \Longleftrightarrow &
v=\left(I+\frac{1}{2}\begin{bmatrix}
-1 & \bar{x}_2^T\\
\bar{x}_2 & -\bar{x}_2\bar{x}_2^T
\end{bmatrix}\right)(u+v).
\end{eqnarray*}
The following argument is similar to Case 3. Let
\[M:=I+\frac{1}{2}\begin{bmatrix}
-1 & \bar{x}_2^T\\
\bar{x}_2 & -\bar{x}_2\bar{x}_2^T
\end{bmatrix}.\]
First note that $I-2M=-I+\begin{bmatrix}
1 & -\bar{x}_2^T\\
-\bar{x}_2 & \bar{x}_2\bar{x}_2^T
\end{bmatrix}$ has the eigenvalue $1$ and $-1$ with multiplicity $n-1$. Hence it is invertible. { For any $t_n\downarrow 0$}, take
$d^n=(d^n_1,d^n_2)$ with
\[d^n_1:=\frac{\|x_2+t_nd_2\|-\|x_2\|}{t_n}, \ d^n_2:=d_2, \ {\rm and}\ w^n:=t_n(\hat{x}+t_n\widehat{d^n}).\]
So $d^n\to d$ and $w^n\to 0$. Then
$(x+t_nd^n,t_nw^n)\in \Omega$ with $x+t_nd^n,t_nw^n\in {\rm bd}\K\backslash \{0\}$. Let
$z^n:=x+t_nd^n-t_nw^n$. Then $z^n\to x$ and
\[
\nabla \Pi_{\K}(z^n)=I+\frac{1}{2}\begin{bmatrix}
-1 & \frac{(z^n_2)^T}{\|z^n_2\|}\\
\frac{z^n_2}{\|z^n_2\|} & -I+\frac{z^n_1}{\|z^n_2\|}(I-\frac{z^n_2}{\|z^n_2\|}\frac{(z^n_2)^T}{\|z^n_2\|})\\
\end{bmatrix}\longrightarrow  I+\frac{1}{2}\begin{bmatrix}
-1 & \bar{x}_2^T\\
\bar{x}_2 & -\bar{x}_2\bar{x}_2^T
\end{bmatrix}=M.
\]
Hence $I-2\nabla \Pi_{\K}(x+t_nd^n-t_nw^n)$ is inverse as $n$ is large enough. Let
\[\alpha(n):=\bigg(I-2\nabla \Pi_{\K}(x+t_nd^n-t_nw^n)\bigg)^{-1}\bigg(\nabla\Pi_{\K}(x+t_nd^n-t_nw^n)-M\bigg)(-u-v).\]
Then
 \[
 \bigg(I-2\nabla \Pi_{\K}(x+t_nd^n-t_nw^n)\bigg)\alpha(n)=\bigg(\nabla\Pi_{\K}(x+t_nd^n-t_nw^n)-M\bigg)(-u-v).
 \]
 Hence
 \begin{eqnarray*}
 \alpha(n)&=&(\nabla\Pi_{\K}(x+t_nd^n-t_nw^n)-M)(-u-v)+2\nabla \Pi_{\K}(x+t_nd^n-t_nw^n)\alpha(n)\\
 &=&\nabla\Pi_{\K}(x+t_nd^n-t_nw^n)(-u+\alpha(n)-v+\alpha(n))-M(-u-v).
 \end{eqnarray*}
 This together with $v=M(u+v)$ yields
 \begin{eqnarray*}
 -v+\alpha(n)= D^*\Pi_{\K}(x+t_nd^n-t_nw^n)(-u+\alpha(n)-v+\alpha(n)).
 \end{eqnarray*}
 So $(u-\alpha(n),v-\alpha(n))\in \widehat{N}_{\Omega}(x+t_nd^n,y+t_nw^n)$. Thus
 $(u,v)\in N^i_{\Omega}((x,y);(d,w))$.

 \underline{Subcase 4.2(iii).} $u\in \mathbb{R}_-\hat{x}$ and $\langle v, \hat{x}\rangle \leq 0$. {For any $t_n\downarrow 0$}, let $w^n:=0$ and $d^n=(d^n_1,d^n_2)$ with $d^n_1:=\frac{\|x_2+t_nd_2\|-x_1}{t_n}$
 and $d^n_2:=d_2$. Then $d^n\to d$, $w^n\to w$, and $x+t_nd^n\in {\rm bd}\K\backslash\{0\}$. Let $z^n:= \hat{x}+t_n\widehat{d^n}$ and
 \[
 u^n:=\frac{u_1}{x_1} z^n \ \ {\rm and}\ \ v^n:=v-t_n\frac{\langle \hat{v},  d^n\rangle }{\|z^n\|}\frac{z^n}{\|z^n\|}.
 \]
 Then $v^n\to v$ and
 $u^n\to \frac{u_1}{x_1}\hat{x}=u$, where $\frac{u_1}{x_1}\leq 0$ and $u=\frac{u_1}{x_1}\hat{x}$ is due to $u\in \mathbb{R}_-\hat{x}$.
 Note that
 $u^n\in \mathbb{R}_- z^n$ and
 \[\langle v^n, z^n\rangle =\langle v,z^n \rangle-t_n\langle \hat{v}, d^n \rangle =\langle \hat{v},x+t_nd^n \rangle-t_n\langle \hat{v}, d^n \rangle=\langle v, \hat{x} \rangle\leq 0.\]
 This means $(u^n,v^n)\in \widehat{N}_{\Omega}(\widehat{z^n},0)=\widehat{N}_{\Omega}(x+t_nd^n,y+t_nw^n)$.
 So $(u,v)\in N_{\Omega}^i((x,y);(d,w))$.

 \underline{Subcase 4.3.} $d\perp \hat{x}$ and $w\in \mathbb{R}_+\hat{x}\backslash\{0\}$. In this case, we will show that
\[ N^i_{\Omega}((x,y);(d,w))=N_{\Omega}((x,y);(d,w))= \{(u,v)|\ u\in \mathbb{R} \hat{x}, v\perp \hat{x}\}.\label{eqn1} \]
 Take $(u,v)\in N_{\Omega}((x,y);(d,w))$. Then there exist sequences $t_n\downarrow 0, (d^n,w^n)\to (d,w), (u^n,v^n)\to (u,v)$ such that $(u^n,v^n)\in \widehat{N}_{\Omega}(x+t_nd^n,t_nw^n)$. Since $x\in {\rm bd}\K\backslash\{0\}$ and $w\not =0$, for $n$ sufficiently large, $0\not =x+t_nd^n \in {\cal K}$ and $0\not = t_nw^n$. It follows that $x+t_nd^n,t_nw^n\in {\rm bd}\K
 \backslash\{0\}$.
 % since $x+t_nd^n$ and $t_nw^n$ are both nonzero and $(x+t_nd^n,t_nw^n)\in \Omega$.
 Hence,
by Proposition \ref{formula-normal cone},  \[
 u^n\perp x+t_nd^n, \ \ v^n\perp t_n w^n, \ \ (x+t_nd^n)_1\widehat{u^n}+(t_nw^n)_1v^n\in \mathbb{R} [x+t_nd^n].
 \]
 Taking the limits yields $u \perp x, v \perp w, \hat{u}\in \mathbb{R}x$, which together with $w\in \mathbb{R}_+\hat{x}\backslash\{0\}$ implies that $v\perp \hat{x}, u\in \mathbb{R}\hat{x}$.  Hence
 \[N_{\Omega}((x,y);(d,w))\subseteq \{(u,v)|\ u\in \mathbb{R} \hat{x}, v\perp \hat{x}\}.\]
 %, i.e., $N_{\Omega}((x,y);(d,w))\subseteq \{(u,v)|\ u\in \mathbb{R} \hat{x}, v\perp \hat{x}\}$.

 Conversely, let $u\in \mathbb{R} \hat{x}$ and $v\perp\hat{x}$. Similarly  to Subcase 4.2(ii), we can prove $(u,v) \in N^i_\Omega((x,y);(d,w))$. The only change is to  take  $w^n:=\frac{w_1}{x_1}(\hat{x}+t_n\widehat{d^n})$ instead of $w^n=t_n(\hat{x}+t_n\widehat{d^n})$.  Since $w=\tau \hat{x}$ for some $\tau>0$, we have $w_2=\frac{w_1}{x_1}(-x_2)$, and hence   $w^n=\frac{w_1}{x_1}(\hat{x}+t_n\widehat{d^n})\rightarrow \frac{w_1}{x_1}\hat{x}=w$.

 \noindent {\bf Case 5.} $x=0$ and $y\in {\rm bd}\K/\{0\}$. The result follows by a symmetric analysis of Case 4.

\noindent {\bf Case 6.} If $x=0$ and $y=0$, then $(d,w)\in T_{\Omega}(x,y)=\Omega$. Using Proposition \ref{lem-directional} yields
 \[
 N_{\Omega}^i((x,y);(d,w))=N_{\Omega}((x,y);(d,w))=N_{\Omega}(d,w).
 \]
 \endproof

%Remark. For a general set, even for convex set we have pointed out the set $\bigcup\limits_{0\neq d\in T_{\Gamma}(x)}N(x;d)$ maybe not equal to $N_{\Gamma}(x)$. Here we point out for the special nonconvex set $\Omega$, these two sets coincide for nonzero $(x,y)$.
%
% \begin{cor}
% For $(x,y)\in \Omega$,
% \[
% \bigcup\limits_{0\neq (d,w)\in T_\Omega(x,y)}N((x,y);(d,w))=N_{\Omega}(x,y), \  \ {\rm if}\ \  (x,y)\neq (0,0)
% \]
%  \begin{equation}\label{co-directional norm}
% \bigcup\limits_{0\neq (d,w)\in T_\Omega(x,y)}N((x,y);(d,w))\subsetneqq N_{\Omega}(x,y), \  \ {\rm if}\ \  (x,y)= (0,0)
% \end{equation}
% \end{cor}
%
% Note that (\ref{co-directional norm}) do not hold by equation. In fact
% for $(x,y)=(0,0)$, then
% \[
%  \bigcup\limits_{0\neq (d,w)\in T_\Omega(x,y)}N((x,y);(d,w))=
% \bigcup\limits_{0\neq (d,w)\in \Omega}N(d,w).
% \]
% Take $(u,v)\in N_{\Omega}(x,y)$. Then $u=0,v\in \mathbb{R}^m$ corresponding to $(d,w)\in {\rm int}\K\times \{0\}$; $u\in \mathbb{R}^m,v=0$ corresponding to $(d,w)\in \{0\}\times {\rm int}\K$; $u\in \mathbb{R}_-\xi, v\in \xi^\circ$ corresponding to $(d,w)=(\hat{\xi},0)$; $u\in \xi^\circ, u\in \mathbb{R}_-\xi$ corresponding to $(d,w)=(0,\hat{\xi})$;
% $u\perp \xi, v\perp \hat{\xi}, \alpha u+(1-\alpha)\hat{v}\in \mathbb{R}\hat{\xi}$ corresponds to $(d,w)=(\alpha \xi, (1-\alpha)\hat{\xi})$. While for $u\in -\K, v\in -\K$, we cannot find nonzero $(d,w)\in \Omega$ satisfying $(u,v)\in N_{\Omega}(d,w)$. For example, take $u=v=(-1,0)$, then $(u,v)\notin \bigcup\limits_{0\neq (d,w)\in \Omega}N_{\Omega}(d,w)$.

\section{Sufficient conditions for error bounds of the second-order cone complementarity system}

In this section, we give verifiable sufficient conditions for the  error bound property of the second-order cone complementarity system (\ref{system-1})-(\ref{system-2}).
First by applying Theorem \ref{Thm4.1} we have the following sufficient conditions based on the  limiting normal cones.
\begin{thm}\label{finalnew}
Given a  point $z^*\in {\cal F}$.  The system (\ref{system-1})-(\ref{system-2}) has a local error bound at $z^*$ if the
NNAMCQ holds at $z^*$:
\[
     \left.
     \begin{array}{l}
    \nabla F(z^*)^T \lambda^F+\displaystyle \sum_{i=1}^J \left \{ \nabla G_i(z^*)^T\lambda^G_i+\nabla H_i(z^*)^T\lambda^H_i\right \}=0, \\
    \lambda^F\in N_\Lambda(F(z^*)), \ \
     (\lambda^G,\lambda^H)\in N_{\Omega}(G(z^*),H(z^*))
     \end{array}
     \right\}\Longrightarrow (\lambda^F,\lambda^G,\lambda^H)=0.
     \]
    Here the exact expression for the limiting  normal cone of $\Omega$ can be found in Proposition \ref{formula-normal cone}.
\end{thm}

%Since the set-valued mapping $M(z):=F(z)-\Lambda$ is metrically subregular at $(z^*,0)$ by \cite{int:loffe}, the normal cone can be estimated by
%$$N_{\cal F} (z^*)\subseteq \{\nabla F(z^*)^T N_\Lambda(F(z^*)).$$ This leads to (ii).
%\endproof

By applying Theorems \ref{ThSuffCondMS} and \ref{PropEquGrSubReg} respectively,
%We now conclude with
we obtain the
  sufficient conditions in Theorems \ref{finaln} and \ref{final}  based on directional limiting normal cone immediately. According to the relationship between the limiting normal cone and the directional limiting normal cone $N_{C}(z;d) \subseteq N_{C}(z)$,
  the sufficient condition based on the directional limiting normal cone is in general weaker than the one based on the limiting normal cone given in Theorem \ref{finalnew}. In fact Example \ref{ex4.1} shows that it is possible that the NNAMCQ does not hold while the sufficient condition in terms of the directional limiting normal cone holds.
  Note that in the following theorem, the formula of the tangent cone and the directional normal cone for a second-order cone complementarity set can be found in Theorems \ref{formula-regular normal cone} and  \ref{formula-directional normal}, respectively. Moreover, the equivalence of the two conditions are due to the directional regularity of the second-order cone complementarity set proved in  Theorem  \ref{formula-directional normal}.
  \begin{thm}\label{finaln}
Given a  point $z^*\in {\cal F}$. Suppose that   for every $0\neq w\in \mathbb{R}^n$ with $\nabla F(z^*)w\in T_{\Lambda}(F(z^*))$, $(\nabla G_i(z^*)w,\nabla H_i(z^*)w)\in T_{\Omega_i}(G_i(z^*),H_i(z^*)), i=1,\dots, J$, one has
     \[
     \left.
     \begin{array}{l}
    \nabla F(z^*)^T \lambda^F+\displaystyle \sum_{i=1}^J \left \{ \nabla G_i(z^*)^T\lambda^G_i+\nabla H_i(z^*)^T\lambda^H_i\right \}=0, \\
    \lambda^F\in N_\Lambda(F(z^*); \nabla F(z^*)w), \\
     (\lambda^G_i,\lambda^H_i)\in N_{\Omega_i}\bigg((G_i(z^*),H_i(z^*));(\nabla G_i(z^*)w,\nabla H_i(z^*)w) \bigg), i=1,\dots, J
     \end{array}
     \right\}\Longrightarrow (\lambda^F,\lambda^G,\lambda^H)=0,
     \]
     or equivalently
    for   every $0\neq w\in \mathbb{R}^n$ with $\nabla F(z^*)w\in T_{\Lambda}(F(z^*))$, $(\nabla G(z^*)w,\nabla H(z^*)w)\in T_{\Omega}(G(z^*),H(z^*))$ one has
     \[
     \left.
     \begin{array}{l}
    \nabla F(z^*)^T \lambda^F+ \nabla G(z^*)^T\lambda^G+\nabla H(z^*)^T\lambda^H=0, \\
    \lambda^F\in N_\Lambda(F(z^*); \nabla F(z^*)w), \\
     (\lambda^G,\lambda^H)\in N_{\Omega}\bigg((G(z^*),H(z^*));(\nabla G(z^*)w,\nabla H(z^*)w) \bigg)
     \end{array}
     \right\}\Longrightarrow (\lambda^F,\lambda^G,\lambda^H)=0.
     \]
%     or equivalently
%     \[
%     \left.
%     \begin{array}{c}
%     \nabla G(\bar{z})^T\lambda^G+\nabla H(\bar{z})^T\lambda^H=0, \\
%     (\lambda^G,\lambda^H)\in \bigcup\limits_{0\neq u\in T^{\rm lin}(z^*)} N_{\Omega}\bigg((G(z^*),H(z^*));(\nabla G(z^*)u,\nabla H(z^*)u) \bigg)
%     \end{array}
%     \right\}\Longrightarrow (\lambda^G,\lambda^H)=0.
%     \]
%     where $T^{\rm lin}(z^*):=\{u|\, (\nabla G(z^*)u,\nabla H(z^*)u)\in T_{\Omega}(G(z^*),H(z^*))\}$.
Then the system (\ref{system-1})-(\ref{system-2}) has a local error bound at $z^*$.

\end{thm}

\begin{thm}\label{final}
Given a  point $z^*\in {\cal F}$. Suppose that the set-valued mapping $M(z):=F(z)-\Lambda$ is metrically subregular at $(z^*,0)$. Further assume that  for every $0\neq w\in \mathbb{R}^n$ with $\nabla F(z^*)w\in T_{\Lambda}(F(z^*))$, $(\nabla G_i(z^*)w,\nabla H_i(z^*)w)\in T_{\Omega_i}(G_i(z^*),H_i(z^*)), i=1,\dots, J$, one has
     \[
     \left.
     \begin{array}{l}
    \nabla F(z^*)^T \lambda^F+\displaystyle \sum_{i=1}^J \left \{ \nabla G_i(z^*)^T\lambda^G_i+\nabla H_i(z^*)^T\lambda^H_i\right \}=0, \\
    \lambda^F\in N_\Lambda(F(z^*); \nabla F(z^*)w), \\
     (\lambda^G_i,\lambda^H_i)\in N_{\Omega_i}\bigg((G_i(z^*),H_i(z^*));(\nabla G_i(z^*)w,\nabla H_i(z^*)w) \bigg), i=1,\dots, J
     \end{array}
     \right\}\Longrightarrow (\lambda^G,\lambda^H)=0,
     \]
     or equivalently
    for   every $0\neq w\in \mathbb{R}^n$ with $\nabla F(z^*)w\in T_{\Lambda}(F(z^*))$, $(\nabla G(z^*)w,\nabla H(z^*)w)\in T_{\Omega}(G(z^*),H(z^*))$ one has
     \[
     \left.
     \begin{array}{l}
    \nabla F(z^*)^T \lambda^F+ \nabla G(z^*)^T\lambda^G+\nabla H(z^*)^T\lambda^H=0, \\
    \lambda^F\in N_\Lambda(F(z^*); \nabla F(z^*)w), \\
     (\lambda^G,\lambda^H)\in N_{\Omega}\bigg((G(z^*),H(z^*));(\nabla G(z^*)w,\nabla H(z^*)w) \bigg)
     \end{array}
     \right\}\Longrightarrow (\lambda^G,\lambda^H)=0.
     \]
%     or equivalently
%     \[
%     \left.
%     \begin{array}{c}
%     \nabla G(\bar{z})^T\lambda^G+\nabla H(\bar{z})^T\lambda^H=0, \\
%     (\lambda^G,\lambda^H)\in \bigcup\limits_{0\neq u\in T^{\rm lin}(z^*)} N_{\Omega}\bigg((G(z^*),H(z^*));(\nabla G(z^*)u,\nabla H(z^*)u) \bigg)
%     \end{array}
%     \right\}\Longrightarrow (\lambda^G,\lambda^H)=0.
%     \]
%     where $T^{\rm lin}(z^*):=\{u|\, (\nabla G(z^*)u,\nabla H(z^*)u)\in T_{\Omega}(G(z^*),H(z^*))\}$.
Then the system (\ref{system-1})-(\ref{system-2}) has a local error bound at $z^*$.

\end{thm}

 In order to use  Theorem \ref{final}, the set-valued mapping $M(z):=F(z)-\Lambda$ should satisfy the metric subregularity.  For convenience, we summarize some prominent sufficient conditions for the case of an equality and inequality system in the following theorem.
  It is well known that in  Theorem \ref{Thm7.2} { (ii)$\Longrightarrow$(iii)$\Longleftrightarrow$(iv)$\Longrightarrow$ (v) or (vi), (i)$\Longrightarrow$ (v) and  (i)$\Longrightarrow$ (vi).}

\begin{thm}[Sufficient conditions for MS for the equality and inequality system] \label{Thm7.2} Let $z^*$ be a feasible point to the system $g(z)\leq 0, h(z)=0$, where $g:\mathbb{R}^n\to \mathbb{R}^p$, $h:\mathbb{R}^n\to \mathbb{R}^q$ are differentiable.  Then the set-valued mapping $M(z):=(g(z),h(z))- {\mathbb{R}_-^p}\times \{0\}^q$ is metrically subregular at $(z^*,0)$ under one of the following conditions.
\begin{itemize}
\item[\rm (i)] Linearity constraint qualification (Linear CQ) holds: $h, g$ are affine.
\item[\rm (ii)] {Linear independence constraint qualification (LICQ) holds: $\{\nabla g_i(z^*), \nabla h_j(z^*)| i\in I_g, j=1,\dots,q\}$ are linearly independent.}
\item[\rm (iii)] The Mangasarian-Fromovitz constraint qualification (MFCQ) holds at $z^*$:   $\{\nabla h_i(z^*)|i=1,\dots,q\}$ are linearly independent and there exists $d\in \mathbb{R}^n$ such that
$\nabla h_i(z^*)d=0$ for all $i=1,\dots,q$ and $\nabla g_i(z^*)d<0$ for all $i\in I_g$.
\item[\rm (iv)] NNAMCQ  holds at $z^*$:
\begin{eqnarray}
&& \nabla g(z^*)^T\lambda^g+\nabla h(z^*)^T\lambda^h=0, \ \
\lambda^g\geq 0, \langle \lambda^g, g(z^*)\rangle =0\label{(22)}\\
&&
  \Longrightarrow (\lambda^g,\lambda^h)=0. \nonumber
\end{eqnarray}
\item[\rm (v)] Quasinormality holds at $z^*$ (\cite[Corollary 5.3]{GuoYeZhang}):
\begin{eqnarray*}
&& \left \{\begin{array}{l}
(\ref{(22)})\mbox{ and there exists a sequence }  \{z^k\} \mbox{ converging to } z^* \mbox{ such that for each }  k, \\
\lambda^g_i> 0 \Longrightarrow {g_i(z^k) >0}, \lambda^h_i \not =0 \Longrightarrow \lambda^h_i h_i(z^k) >0
\end{array} \right.  \\
&& \Longrightarrow (\lambda^g,\lambda^h)=0.
\end{eqnarray*}
\item[\rm (vi)] The  relaxed constant positive linear dependence condition (RCPLD) holds at $z^*$ (\cite[Theorem 4.2]{GuoZhangLin}): Let ${\cal{J}}\subseteq \{1,\dots,q\}$ and $\{\nabla h_j(z^*)|j\in {\cal J}\}$ be a basis for ${\rm span}\{\nabla h_j(z^*)|j=1,\dots,q\}$. There exists $\delta>0$ such that
\begin{itemize}
\item $\{\nabla h_j(z)\}_{j=1}^q$ has the same rank for each $z\in B_\delta(z^*)$;
\item   for each ${\cal I}\subseteq I_g$, if $\{\nabla g_i(z^*),\nabla h_j(z^*)|i\in {\cal I},j\in {\cal J}\}$ is positively linear dependent, then $\{\nabla g_i(z),\nabla h_j(z)|i\in {\cal I},j\in {\cal J}\}$ is linear dependent for each $z\in B_\delta(z^*)$.
\end{itemize}
\item[\rm (vii)] There are no nonzero direction in the  linearized cone (\cite[Corollary 1]{GfrKl16}):
$$  \nabla g_i(z^*)d\leq 0,  i\in I_g, \nabla h_i(z^*)d=0,
i=1,\dots,q \Longrightarrow d=0.$$
%\item  First-order sufficient condition for metric subregularity (FOSCMS) (\cite[Corollary 1]{GfrKl16}): For every $0\neq u\in \mathbb{R}^n$ with $\nabla g_i(z^*)u\leq 0$ for $i\in I_g(z^*)$ and $\nabla h_j(z^*)u=0$ for $j=1,\dots,q$, one has
%\begin{eqnarray*}
%&& \nabla g(z^*)^T\lambda^g+\nabla h(z^*)^T\lambda^h=0, \ \
%\lambda^g\geq 0, \langle \lambda^g, g(z^*)\rangle =0, \langle \lambda^g,\nabla g(z^*)u\rangle =0\\
%&&  \Longrightarrow (\lambda^g,\lambda^h)=0
%\end{eqnarray*}
\item[\rm (viii)] Second-order sufficient condition for metric subregularity (SOSCMS) (\cite[Corollary 1]{GfrKl16}): For every $0\neq w\in \mathbb{R}^n$ with $\nabla g_i(z^*)w\leq 0$ for $i\in I_g$ and $\nabla h_i(z^*)w=0$ for $i=1,\dots,q$, one has
\begin{eqnarray*}
\left.\begin{array}{l}
\nabla g(z^*)^T\lambda^g+\nabla h(z^*)^T\lambda^h=0, \ \lambda^g\geq 0, \langle \lambda^g, g(z^*)\rangle =0
\\
 w^T\nabla^2 ((\lambda^g)^T g)(z^*)w+w^T\nabla^2 ((\lambda^h)^Th)(z^*)w\geq 0
\end{array}
\right \}
\Longrightarrow (\lambda^g,\lambda^h)=0.
\end{eqnarray*}
\end{itemize}
\end{thm}

The following example shows that if there exists $0\neq w\in \mathbb{R}^n$ with $\nabla g_i(z^*)w\leq 0$ for $i\in I_g$ and $\nabla h_i(z^*)w=0$ for $i=1,\dots,q$, then SOSCMS is weaker than NNAMCQ, or equivalently MFCQ.
 \begin{example}
 Let $g_1(z)=z_1-z_2^2$, $g_2(z)=z_1^2-z_2$, and $h(z)=z_1$. At $z=(0,0)$, consider
 \[
 \lambda^{g_1}\nabla g_1(z)+\lambda^{g_2}\nabla g_2(z)+\lambda^h \nabla h(z)=0
 \]
 with $\lambda^{g_1},\lambda^{g_2}\geq 0$. Then we can take $(\lambda^{g_1},\lambda^{g_2},\lambda^{h})=(1,0,-1)\not = (0,0,0)$. So MFCQ fails at $z=(0,0)$.
 Let $w$ satisfying $\nabla h(z)w=0$. Then $w_1=0$, and hence from
 \[
 \lambda^{g_1}w^T\nabla^2 g_1(z)w+\lambda^{g_2}w^T\nabla^2 g_2(z)w\geq 0,
 \]
 we have $-2\lambda^{g_1}w_2^2\geq 0$, so $\lambda^{g_1}\leq 0$ and hence $\lambda^{g_1}=0$. Consequently, $\lambda^{g_2}=0$ and $\lambda^h=0$. So SOSCMS holds at $z=(0,0)$.
 \end{example}

Our sufficient condition Theorem \ref{final} provides a sufficient condition for metric subregularity for the very general system (\ref{system-1})-(\ref{system-2}).
There may exist more than one way to split a system and this provides flexibility in using Theorem \ref{final}. For example,
suppose that  a second-order cone complementarity system consists only  (\ref{system-1}). Suppose some of the subsystems, without loss of generality,
$$ (G_i(z),H_i(z))  \in \Omega_i,\ \  i=1,\dots, s,$$
where $s\leq J$
is metric  subregular at $(z^*,0)$. Then one can split the original system
(\ref{system-1}) as
\begin{eqnarray}
&&{\cal K}_i \ni G_i(z)  \perp H_i(z)  \in {\cal K}_i \quad i=s+1,\dots, J, \nonumber\\
&&  F(z) \in \Lambda,\label{vectorc}
\end{eqnarray}
where
$$F(z):= (G_1(z),H_1(z)) \times \dots \times   (G_s(z),H_s(z)) \ \mbox{ and } \ \Lambda:= \Omega_1\times \dots \times \Omega_s,$$
and use Theorem \ref{final}.
In particular, since ${\cal K}_i$ with $m_i=1$ is  the set of  nonnegative reals $\mathbb{R}_+$,
$\Omega_i$ with $m_i=1$ is equal to the vector complementarity cone  ${\Theta}:=\{(x,y)\in \mathbb{R}^2| x\geq 0,
y\geq 0, x^Ty=0\}$. Without loss of generality, assuming $m_i=1$ for $i=1,\dots, s$, then the system $F(z)\in \Lambda$ given in (\ref{vectorc})  is then equal to the vector complementarity system
$$  F(z)\in {\Theta}^s,$$
where ${\Theta}^s:=\{(a,b)\in \mathbb{R}^{2s}| a\geq 0,b\geq 0, a^Tb=0\}$. We now summarize  some prominent sufficient conditions for metric subregularity for the vector complementarity system in the following theorem.

\begin{thm}[Sufficient conditions for MS for a complementarity system] Let $z^*$ be a feasible point to the vector complementarity system $g(z)\leq 0, h(z)=0,0\leq \phi(z)\perp \psi(z) \geq 0$,
where $g:\mathbb{R}^n\to \mathbb{R}^p$, $h:\mathbb{R}^n\to \mathbb{R}^q$, $\phi:\mathbb{R}^n\to \mathbb{R}^s$, $\psi:\mathbb{R}^n\to \mathbb{R}^s$ are continuously differentiable.
%For a function $g$, denote by
%\begin{eqnarray*}
%&& I_g= I_g(z^*):=\{i| g_i(z^*)=0\}
%\end{eqnarray*}
%the active index set.
Then, the set-valued mapping $M(z):=(g(z),h(z), \phi(z), \psi(z))- {\mathbb{R}_-^p}\times \{0\}^q\times {\Theta}^s$ is metrically subregular at $(z^*,0)$ under one of the following conditions:
\begin{itemize}
\item[\rm (i)] Linearity CQ holds: $g, h, \phi, \psi$ are affine.
\item[\rm (ii)]  MPEC LICQ holds:  $\{\nabla g_i(z^*) (i\in I_g), \nabla h_i(z^*) (i=1,\dots,q), \nabla \phi_i(z^*)  (i\in I_\phi), \nabla \psi_i(z^*)  (i\in I_{\psi})\}$ are linearly independent.
\item[\rm (iii)] MPEC NNAMCQ holds at $z^*$:
\begin{eqnarray}
&& \nabla g(z^*)^T\lambda^g+ \nabla h(z^*)^T\lambda^h{+ \nabla \phi(z^*)^T\lambda^\phi+\nabla \psi(z^*)^T\lambda^\psi}=0, \label{(23)}\\
&& \lambda^g \geq 0, \lambda^g_i=0, \forall   i\not \in I_g, \lambda_i^\phi=0, \forall i \not \in I_\phi, \lambda_i^\psi=0, \ \forall i \not \in I_\psi, \label{(24)}\\
&& \mbox{ either } {\lambda^\phi_i <0}, {\lambda^\psi_i<0} \mbox{ or } \lambda^\phi_i \lambda^\psi_i=0 \ \ \forall i\in I_\phi\cap I_\psi \label{(25)}\\
&&  \Longrightarrow (\lambda^g,\lambda^h, \lambda^\phi,\lambda^\psi)=0. \nonumber
\end{eqnarray}
%MPEC MFCQ holds at $z^*$ (\cite{ss}): \begin{eqnarray}
%&& \nabla \phi(z^*)^T\lambda^\phi+\nabla \psi(z^*)^T\lambda^\psi=0 \label{(23)}\\
%&& \lambda_i^\phi=0, \forall i \not \in I_\phi(z^*), \lambda_i^\psi=0 \forall i \not \in I_\psi(z^*) \label{(24)}\\
%&&  \Longrightarrow (\lambda^\phi,\lambda^\psi)=0 \nonumber
%\end{eqnarray}

%\item Quasinormality holds at $z^*$.
%\item The  relaxed constant positive linear dependence condition (RCPLD) holds at $z^*$ %(\cite{GuoZhangLin}): Let ${\cal{J}}\subseteq \{1,\dots,q\}$ and $\{\nabla h_j(z^*)|j\in J\}$ is a %basis for ${\rm span}\{\nabla h_j(z^*)|j=1,\dots,q\}$. There exists $\delta>0$ such that
%\begin{enumerate}
%\item[\rm(i)] $\{\nabla h_j(z)\}_{j=1}^q$ has the same rank for each $z\in B(z^*,\delta)$;
%\item[\rm(ii)] for each ${\cal I}\subseteq I_g(z^*)$, if $\{\nabla g_i(z^*),\nabla h_j(z^*)|i\in %%{\cal I},j\in {\cal J}\}$ is positively linear dependent, then $\{\nabla g_i(z),\nabla h_j(z)|i\in %%{\cal I},j\in {\cal J}\}$ is linear dependent for each $z\in B(z^*,\delta)$.
%\end{enumerate}
\item[\rm (iv)] MPEC quasi-normality holds at $z^*$:
\begin{eqnarray*}
&& \left \{ \begin{array}{l}
 (\ref{(23)})-(\ref{(25)}) \mbox{ and there exists a sequence }  \{z^k\} \mbox{ converging to } z^* \mbox{ such that for each }  k, \\
\lambda^g_i>0 \Rightarrow \lambda^g_i g_i(z^k) >0, \lambda^h_i\not =0 \Rightarrow \lambda^h_i h_i(z^k)
>0,\\
 \lambda^\phi_i\not =0 \Rightarrow \lambda^\phi_i \phi_i(z^k) <0, \lambda^\psi_i \not =0 \Rightarrow \lambda^\psi_i \psi_i(z^k) <0,  \end{array}
\right. \\
&&  \Longrightarrow (\lambda^g,\lambda^h,\lambda^\phi,\lambda^\psi)=0. \nonumber
\end{eqnarray*}
\item[\rm (v)] There is no nonzero direction in the MPEC  linearized cone: $ \mathcal{L}^{MPEC}(z^*)=\{0\}$ where
\begin{eqnarray*}
  \mathcal{L}^{MPEC}(z^*)&:=&\left \{w \left | \begin{array}{l} \nabla g_i(z^*)w\leq 0, i\in I_g, \nabla h_i(z^*)w=0,i=1,\dots, q\\
  \nabla  \phi_i(z^*)w=0, i \in I_\psi^c, \nabla \psi_i(z^*)w=0, i \in I_\phi^c\\
0\leq \nabla \phi_i(z^*)w\perp \nabla \psi_i(z^*) w\geq 0 \quad i\in I_\phi\cap I_\psi
\end{array}
\right . \right \}.
\end{eqnarray*}
\item[\rm (vi)]  The set-valued mapping  $M_1(z):=(g(z),h(z))- {\mathbb{R}_-^p}\times \{0\}^q$  is metrically subregular at $(z^*,0)$ and for every
$0\not =w\in  \mathcal{L}^{MPEC}(z^*)$,  one has
\begin{eqnarray}
&& (\ref{(23)})-(\ref{(24)}) \mbox{ and }
% \sum_{i\in I_\phi}  \lambda_i^\phi \nabla \phi_i(z^*) u+\sum_{i\in  I_\psi} \lambda_i^\psi \nabla \psi_i(z^*) u=0
(\lambda_i^\phi, \lambda_i^\psi) \in N_{\Theta}(\nabla \phi_i(z^*)w,\nabla \psi_i(z^*)w) \ \forall i\in I_\phi\cap I_\psi \label{FOSCMS} \\
&&  \Longrightarrow (\lambda^\phi,\lambda^\psi)=0,\nonumber
\end{eqnarray}
where $N_{\Theta}(x, y)$ is given in (\ref{normalc}).
\item[\rm (vii)] The set-valued mapping
 $M_1(z):=(g(z),h(z))- { \mathbb{R}_-^p}\times \{0\}^q$  is metrically subregular at $(z^*,0)$ and for every
$0\not =w\in  \mathcal{L}^{MPEC}(z^*)$,
 one has
\[
 (\ref{FOSCMS})
\mbox{ and } w^T\nabla_z^2 L^0(z^*,\lambda^g,\lambda^h,\lambda^\phi,\lambda^\psi)w \geq 0  \Longrightarrow (\lambda^\phi,\lambda^\psi)=0,
\]
where
$${L^0}(z, \lambda^g,\lambda^h,\lambda^\phi,\lambda^\psi):=g(z)^T\lambda^g+h(z)^T\lambda^h+\phi(z)^T\lambda^\phi+\psi(z)^T\lambda^\psi.$$
\end{itemize}
\end{thm}
\beginproof (i) follows from the corollary in  \cite[page 210]{Robinson81}. (ii) is stronger than (iii), which is further stronger than (iv).
(iv) follows from \cite[Theorem 5.2]{GuoYeZhang}.
(v) is the trivial case of (vi).
(vi) and (vii)  follow from \cite[Theorem 2.6]{Gfr14} and the well-known fact that the limiting normal cone of the complementarity cone ${\Theta}$ is
equal to
\begin{equation}
N_{\Theta}(x, y)=\left \{(u,v)\in \mathbb{R}^2 \left |\begin{array}{ll}
u =0 & \mbox{ if } x>0 \\
v =0 & \mbox{ if } y>0 \\
 \mbox{either }u<0, v<0 \mbox{ or } uv=0 &\mbox{ if } x=y=0\end{array}\right .
\right  \}.\label{normalc}
\end{equation}

\endproof

We now consider the following SOCCP
\begin{eqnarray}
&&{\cal K}  \ni G(z)  \perp H(z)  \in {\cal K},\label{SOC} \\
&&  g(z)\leq 0, h(z)=0,0\leq \phi(z)\perp \psi(z) \geq 0, \label{CS}
\end{eqnarray}
 where  the second-order cone complementarity system (\ref{SOC}) is defined as in (\ref{system-1}) and  $g:\mathbb{R}^n\to \mathbb{R}^p$, $h:\mathbb{R}^n\to \mathbb{R}^q$, $\phi:\mathbb{R}^n\to \mathbb{R}^s$, $\psi:\mathbb{R}^n\to \mathbb{R}^s$,
 $G:\mathbb{R}^n\to \mathbb{R}^{m}$, $H:\mathbb{R}^n\to \mathbb{R}^{m}$ are continuously differentiable.
 Let the linearized cone of the system (\ref{SOC})-(\ref{CS}) be
 \begin{eqnarray*}
  {\mathcal{L}(z^*)} &:=& \left\{ w \left | \begin{array}{l} \nabla g_i(z^*)w\leq 0, i\in I_g, \nabla h_i(z^*)w=0,i=1,\dots, q\\
  \nabla  \phi_i(z^*)w=0, i \in I_\psi^c, \nabla \psi_i(z^*)w=0, i \in I_\phi^c\\
0\leq \nabla \phi_i(z^*)w\perp \nabla \psi_i(z^*) w\geq 0 \quad i\in I_\phi\cap I_\psi \\
(\nabla G_i(z^*)w,\nabla H_i(z^*)w)\in T_{\Omega_i}(G_i(z^*),H_i(z^*)), i=1,\dots, J
\end{array}
\right. \right \},
\end{eqnarray*}
where $T_{\Omega_i}$ can be calculated as in  Theorem \ref{formula-regular normal cone}. Based on the results we obtained, we now derive a  sufficient condition for error bounds for the system (\ref{SOC})-(\ref{CS}) that are explicitly verifiable based on the initial data.

\begin{thm}
Given a  point $z^*\in {\cal F}$. Suppose that the complementarity system (\ref{CS}) is   metrically subregular at $(z^*,0)$. Further assume that either $ \mathcal{L}(z^*)=\{0\}$, or $ \mathcal{L}(z^*)\not =\{0\}$ and for every $0\neq w\in \mathcal{L}(z^*)$ one has
     \begin{eqnarray*}
   &&  \left \{
     \begin{array}{l}
\displaystyle  \nabla g(z^*)^T\lambda^g+ \nabla h(z^*)^T\lambda^h {+ \nabla \phi(z^*)^T\lambda^\phi+\nabla \psi(z^*)^T\lambda^\psi + \nabla G(z^*)^T\lambda^G+\nabla H(z^*)^T\lambda^H=0}, \\
\lambda^g \geq 0, \lambda^g_i=0 \  \forall  \ i\not \in I_g, \lambda_i^\phi=0 \  \forall i \not \in I_\phi, \lambda_i^\psi=0 \ \forall i \not \in I_\psi, \\
%\mbox{ either } \lambda^\phi_i>0, \lambda^\psi_i>0 \mbox{ or } \lambda^\phi_i \lambda^\psi_i=0 \ \ \forall i\in I_\phi\cap I_\psi,\\
(\lambda_i^\phi, \lambda_i^\psi) \in N_{\Theta}(\nabla \phi_i(z^*)w,\nabla \psi_i(z^*)w) \ \forall i\in I_\phi\cap I_\psi,\\
%\sum_{i\in I_\phi}  \lambda_i^\phi \nabla \phi_i(z^*) u+\sum_{i\in  I_\psi} \lambda_i^\psi \nabla \psi_i(z^*) u=0,  \
     (\lambda^G_i,\lambda^H_i)\in N_{\Omega_i}\bigg((G_i(z^*),H_i(z^*));(\nabla G_i(z^*)w,\nabla H_i(z^*)w) \bigg), \ \forall  i=1,\dots, J,
     \end{array}
     \right .\\
     && \Longrightarrow (\lambda^G,\lambda^H)=0,
     \end{eqnarray*}
     where $N_{\Theta}(\cdot)$ is given in (\ref{normalc}) and $N_{\Omega_i}((x,y);(d,w))$  can be calculated as in Theorem  \ref{formula-directional normal}.
Then the system (\ref{SOC})-(\ref{CS}) has a local error bound at $z^*$. That is,
there exist a constant $\kappa >0$ and $\delta>0$ such that
 \begin{eqnarray*}
  d(z, {\cal F})\leq \kappa \left \{\|h(z)\|+\|g_+(z)\|+\sum_{i=1}^s d_{\Theta}(\phi_i(z),\psi_i(z)) +\sum_{i=1}^J d_{\Omega_i} (G_i(z),H_i(z))\right \},
  \forall z\in B_\delta(z^*).
 \end{eqnarray*}
 %where ${\cal C}:=\{(a,b)\in \mathbb{R}^2| a\geq 0,b\geq 0, ab=0\}$.
\end{thm}
\beginproof
To prove the result, we take $F(z):=(g(z),h(z), \phi(z), \psi(z))$ and $\Lambda:= {\mathbb{R}_-^p}\times \{0\}^q\times {\Theta}^s$ and apply  Theorem \ref{final}. Since the sets ${\mathbb{R}_-^p}$, $\{0\}^q$ are convex and ${\Theta}^s$ is directionally regular, we have
\begin{eqnarray*}
&&T_\Lambda(F(z^*))= T_{\mathbb{R}_-^p}(g(z^*))\times T_{ \{0\}^q}(h(z^*))\times T_{{\Theta}^s}(\phi(z^*),\psi(z^*)),\\
&&N_\Lambda(F(z^*); \nabla F(z^*)w)\\
&&= {N_{\mathbb{R}_-^p}(g(z^*);\nabla g(z^*)w)}\times N_{ \{0\}^q}(h(z^*);\nabla h(z^*)w)\times N_{{\Theta}^s}((\phi(z^*),\psi(z^*)); (\nabla \phi(z^*)w, \nabla \phi(z^*)w)),\\
&&N_{{\Theta}^s}\big((\phi(z^*),\psi(z^*)); (\nabla \phi(z^*)w, \nabla \phi(z^*)w)\big)=\Pi_{i=1}^s  N_{{\Theta}}((\phi_i(z^*),\psi_i(z^*)); (\nabla \phi_i(z^*)w, \nabla \phi_i(z^*)w)).
\end{eqnarray*}
Moreover by Proposition \ref{lem-directional} and Theorem \ref{formula-directional normal}, for all $w $ such that
$(\nabla \phi_i(z^*) w, \nabla \psi_i (z^*)w)\in  T_{{\Theta}}(\phi_i(z^*),\psi_i(z^*))$, we have
\begin{eqnarray*}
 \lefteqn{ N_{{\Theta}}((\phi_i(z^*),\psi_i(z^*)); (\nabla \phi_i(z^*) w, \nabla \phi_i (z^*)w)) }\\
 && =\left\{  \begin{array}{ll}
\mathbb{R}\times \{0\}  & \mbox{  if } \phi_i(z^*)=0,\psi_i(z^*)>0,\\
  \{0\} \times \mathbb{R}  & \mbox{  if } \phi_i(z^*)>0,\psi_i(z^*)=0,\\
   N_{{\Theta}}(\nabla \phi_i(z^*) w, \nabla \psi_i (z^*)w) & \mbox{  if } \phi_i(z^*)=0,\psi_i(z^*)=0.
  \end{array}
  \right.
\end{eqnarray*}
By using the tangent cone formula in Corollary \ref{co-tangent-cone} and the limiting normal cone formula in (\ref{normalc}), we derive the result.
\endproof

%since the second-order cone with dimension less or equal to two is a convex polyhedra and hence the corresponding second-order cone complementarity set is a union of finitely many convex polyhedral sets. In this case, we may put the system ${\cal K} _i \ni G_i(z)  \perp H_i(z)  \in {\cal K}_i$ where $m_i\leq 2$ into the part $F(z)\in \Lambda$ since we may be able to use the polyhedral structure in verifying the metric subregularity.

\noindent {\bf Acknowledgments.} {The authors are grateful to the two anonymous referees for their helpful comments and suggestions.}

\end{document}